\documentclass{article}
\usepackage{arxiv}
\usepackage[utf8]{inputenc} % allow utf-8 input
\usepackage[T1]{fontenc}    % use 8-bit T1 fonts
\usepackage{hyperref}       % hyperlinks
\hypersetup{
	pdftitle={Title},
	pdfsubject={Computational Mechanics},
	pdfauthor={L.~Maurer},
	pdfkeywords={Computational Mechanics},
	colorlinks=false,
	hidelinks,
}
\usepackage{url}            % simple URL typesetting
\usepackage{booktabs}       % professional-quality tables
\usepackage{longtable}
\usepackage{colortbl}
\usepackage{overpic}
\usepackage{amsfonts}       % blackboard math symbols
\usepackage{amsmath}
\usepackage{amssymb}
\usepackage{mathtools}
\usepackage{cancel}
\usepackage{xfrac}			% compact symbols for 1/2, etc.
\usepackage{microtype}      % microtypography
\usepackage{cleveref}       % smart cross-referencing
\usepackage{graphicx}
\usepackage{doi}
\usepackage{subfig}
\usepackage{placeins}
\usepackage{float}
\usepackage{scalerel}
\usepackage{diagbox}		% diagonal in table cell
\usepackage{arydshln}		% dashed lines in table
\usepackage{multirow}
\usepackage[bottom,symbol*]{footmisc}	% Footnotes below floats
\usepackage[noadjust]{cite}
\usepackage{comment}
\usepackage{bbold}
\usepackage{makecell}

\usepackage{xcolor}

\definecolor{OvGU-MB}{cmyk/RGB/HTML}{0,1,0,0.6/122,0,60/7A003F}
\definecolor{Col_red}{HTML}{F09EA7}
\definecolor{Col_ora}{HTML}{F6CA94}
\definecolor{Col_yel}{HTML}{FAFABE}
\definecolor{Col_gre}{HTML}{C1EBC0}
\definecolor{Col_blu}{HTML}{C7CAFF}
\definecolor{Col_vio}{HTML}{CDABEB}
\definecolor{Col_pin}{HTML}{F6C2F3}
\definecolor{Col_gra}{HTML}{F2F2F2}
\newcommand{\tovgu}[1]{\textcolor{OvGU-MB}{#1}}

\usepackage{pgfplots}
\pgfplotsset{compat=1.17}
\usepackage{tikz}
\usepackage{tikzscale}
\newlength\tkwi
\newlength\tkhe
\usetikzlibrary{external}
\newcommand{\smarttikz}[1]{%
  \IfFileExists{figs/PDF-#1.pdf}{%
    \includegraphics{figs/PDF-#1.pdf}%
  }{%
    \tikzsetnextfilename{PDF-#1}%
    \input{tikz/#1.tikz}%
  }%
}

\usepackage{listings}
\usepackage{matlab-prettifier}
\usepackage{caption}
\usepackage{subcaption}

\lstdefinestyle{fortranstyle}{
  language=Fortran,
  basicstyle=\ttfamily\scriptsize,
  keywordstyle=\bfseries,
  commentstyle=\itshape\color{gray},
  numbers=right,
  numberstyle=\tiny\color{black},
  stepnumber=1,
  numbersep=8pt,
  frame=single,
  framerule=0.4pt,
  rulecolor=\color{black!30},
  framesep=5pt,
  xleftmargin=0pt,
  xrightmargin=2.8em,
  framexleftmargin=0pt,
  framexrightmargin=2.8em,
  breaklines=true,
  breakatwhitespace=false,
  columns=fullflexible,
  keepspaces=true,
  showstringspaces=false,
  tabsize=2,
  captionpos=t
}
\let\originalleft\left
\let\originalright\right
\def\left#1{\mathopen{}\originalleft#1}
\def\right#1{\originalright#1\mathclose{}}
\graphicspath{{./figs/}}
\title{Implementation of Hyperelastic Physics-Augmented Neural Networks in the Explicit Finite Element Codes \textsc{Simcenter Radioss} and \textsc{OpenRadioss} with Applications to Impact Events}
\author{
	\href{https://orcid.org/0000-0002-5490-7981}{\includegraphics[scale=0.06]{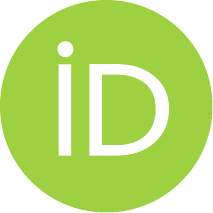}\hspace{1mm}L.~Maurer}\thanks{Corresponding author: \texttt{\href{mailto:lukas.maurer@ovgu.de}{lukas.maurer@ovgu.de}}}\;,
    \href{https://orcid.org/0000-0001-8774-9732}{\includegraphics[scale=0.06]{./logo/orcid.pdf}\hspace{1mm}S.~Eisentr\"ager} \& 
    \href{https://orcid.org/0000-0001-8997-3818}{\includegraphics[scale=0.06]{./logo/orcid.pdf}\hspace{1mm}D.~Juhre}\\
	Institute of Materials, Technologies and Mechanics\\
	Otto von Guericke University Magdeburg\\
	Magdeburg, D-39106
\And
    \href{https://orcid.org/0000-0003-1683-6761}{\includegraphics[scale=0.06]{./logo/orcid.pdf}\hspace{1mm}M.~Bulla}\\
    Siemens Digital Industries\\
    Am Kabellager 9\\
    K\"oln, D-51063
}
\renewcommand{\shorttitle}{PANNs for Explicit Dynamics in \textsc{OpenRadioss}}
\newcommand{\identMat}{\ensuremath{\overset{\raisebox{-0.375ex}{\scalebox{0.5}{<2>}}}{\mathbf{I}}}}
\newcommand{\revI}[1]{{\color{black}#1}}

\newcommand{\revIV}[1]{{\color{black}#1}}

\usepackage{lineno}
\begin{document}
%
%%%%%%%%%%%%%%%%%%%%%%%%%%\linenumbers
%
\maketitle
\begin{abstract}

Data-driven material modeling techniques have gained significant attention due to their ability to capture complex constitutive behaviors beyond the limitations of classical material models. Physics-augmented neural networks (PANNs), which embed physical constraints directly into their architecture, combine the flexibility of machine learning with the reliability required for engineering simulations. As a result, they have become increasingly relevant in computational material modeling.

This work presents an approach to integrate such network architectures into the explicit finite element solvers \textsc{Simcenter Radioss} and \textsc{OpenRadioss} (Siemens). To this end, a framework for transferring pretrained network architectures and their corresponding parameters to a standalone user material routine is developed. The networks are trained using the \textsc{PyTorch} framework, but the proposed procedure can be adapted easily to other machine-learning frameworks, such as \textsc{TensorFlow}. The key advantage of this approach is that it enables the use of PANNs within the finite element technology available in \textsc{Simcenter Radioss} and \textsc{OpenRadioss}, without requiring the development of specialized solvers.

In addition to the integration strategy, particular emphasis is placed on computational efficiency. The influence of the network architecture on the overall simulation performance is investigated, and strategies for reducing evaluation costs while preserving accuracy in the material response are discussed. In this context, we demonstrate how replacing the commonly used SoftPlus activation function by the SQuarePlus function can reduce the computational cost of the PANN material routine. To facilitate practical adoption, a publicly available {\tovgu{GitHub repository}}\footnote[2]{\href{https://github.com/OVGU-CoMe/PANN_Radioss}{https://github.com/OVGU-CoMe/PANN\_Radioss}} is provided that automates the generation of user material routines in Fortran. The proposed workflow allows users to integrate their own neural-network models with limited manual effort, requiring merely the specification of the network architecture and the corresponding trained parameters.

An example impact simulation demonstrates that the generated PANN user material can reproduce the nonlinear behavior characteristic of hyperelastic materials under large strains. Overall, the presented framework provides a practical route toward the use of machine-learning-based constitutive models in explicit finite element simulations. Moreover, it highlights how both commercial and open-source solvers, such as \textsc{Simcenter Radioss} and \textsc{OpenRadioss}, can serve as a testbed for next-generation physics-based neural networks.
\end{abstract}
\section{Introduction}
\label{sec:Introduction}

For decades, scientists and engineers have sought numerical descriptions of real-world materials. Material modeling aims to represent material behavior mathematically and to compare numerical predictions with experimental observations. The more knowledge is available about the micromechanical behavior of a material, the more physical mechanisms can be incorporated into a constitutive description. In most engineering applications, however, not all physical mechanisms need to be explicitly included. A central task of material modeling is therefore to simplify material behavior in a controlled manner, reducing computational effort while still enabling reliable predictions of stresses and strains.

In this work, we focus on isotropic hyperelastic materials, for which the stored energy can be described by a strain-energy density function $\psi$ \cite{Eisentrager.2025, Eisentrager.2025b}. This energy depends only on the deformation state. Additional phenomena such as viscoelasticity, plasticity, damage, or temperature dependence, which are particularly relevant for rubber-like materials, can be added to such a formulation if required. Hyperelastic constitutive laws may be formulated with varying degrees of physical motivation. Some models transfer micromechanical assumptions to the macroscopic scale, as in the Arruda--Boyce model \cite{Arruda.1993}. In contrast, phenomenological models are mainly designed to reproduce experimental observations, such as the Carroll model \cite{Carroll.2011}.

Neural networks have attracted increasing interest in material modeling in recent years. While early approaches often treated the constitutive relation as a black-box mapping directly from deformation measures to stress, more promising results have been obtained with physics-augmented neural networks (PANNs)\cite{Klein.2022, Kalina.2023, Linden.2023, Linka.2023, Schommartz.2025, Damma.2025b}. In these models, key principles of classical constitutive modeling are embedded directly into the network architecture. A central idea is to use the network to represent the scalar potential $\psi$ instead of predicting stresses directly. Stresses, and if required consistent tangent stiffness matrices, can then be obtained by differentiating this potential, which can be readily achieved using automatic differentiation approaches embedded in all machine learning frameworks. This construction helps to ensure thermodynamic consistency. Additional constraints, such as monotonicity or convexity \cite{Klein.2026}, can be enforced through the choice of network structure, weights, and activation functions. In this way, the flexibility of neural networks can be combined with established requirements from continuum mechanics. Since these physical principles are imposed independently of the training data, PANNs provide a robust basis for finite element applications.

While the theoretical formulation and training of physics-augmented neural networks for hyperelastic materials have been studied extensively, their use in established finite element solvers remains a practical challenge. Most neural-network models are trained in machine-learning frameworks such as PyTorch or TensorFlow, whereas industrial finite element solvers typically require compiled user material routines written in languages such as Fortran or C/C++. Directly coupling a finite element solver to a machine-learning framework is possible in principle, but such approaches introduce additional dependencies, may reduce computational efficiency, and are often difficult to integrate into standard simulation workflows.

The present contribution investigates the integration of PANN-based hyperelastic material models into the explicit finite element solver \textsc{OpenRadioss}. \textsc{OpenRadioss} is the open-source version of Siemens \textsc{Simcenter Radioss}. Since both \textsc{OpenRadioss} and \textsc{Simcenter Radioss} support user material routines, the implementation described here can be used in both solvers. For readability, we refer to \textsc{Radioss} in the following, unless a distinction is required. With this integration, users are provided with a simple way to test and embed their own networks in a state-of-the-art finite element solver environment. The Python code that generates the material model directly in Fortran is provided in a public {\tovgu{GitHub repository}}\footnote[1]{\href{https://github.com/OVGU-CoMe/PANN_Radioss}{https://github.com/OVGU-CoMe/PANN\_Radioss}}. We demonstrate how the network architecture and all trained parameters can be exported to Fortran. The resulting workflow does not require PyTorch, TensorFlow, or any other machine-learning environment during the finite element simulation. As a benchmark problem, we consider an impact simulation involving a small dragon-shaped structure. The results of a classical Carroll constitutive law \cite{Carroll.2011} are compared with two PANN architectures of different complexity. We demonstrate that, for the explicit simulations considered here, the generated network-based material routines produce similar results and comparable computational performance.

In Section~\ref{sec:Simsetup}, we briefly introduce the fundamentals of the explicit simulation approach, the benchmark problem, and the hardware used for the numerical studies. Section~\ref{sec:Matmodel} presents the basic continuum-mechanical equations for hyperelastic material modeling, as well as the structure of the PANN approach and the Carroll model, which serves as a classical constitutive law. Both material descriptions are fitted to the Treloar data~\cite{Treloar.1944}. The implementation strategy is described in Section~\ref{sec:Implementation}. In addition to several listings illustrating the behavior of the material routine, we provide ready-to-use Python code that allows users to generate user material subroutines from their own PANN model definitions. In this context, we also investigate replacing the SoftPlus function $\mathcal{SP}\left(x\right)$, which is commonly used in PANN architectures, with the SQuarePlus function $\mathcal{SQP}\left(x\right)$~\cite{Barron.2021}. The SQuarePlus function satisfies the required smoothness, monotonicity, and convexity properties, while being less expensive to evaluate in the generated material routine. To compare both material models, the benchmark test is run multiple times and the corresponding results are evaluated in Section~\ref{sec:Results}. Since the trained PANN model is transferred from Python to Fortran, the simulation runtimes of all material formulations can be compared consistently in Section~\ref{sec:Runtime}. Finally, the main findings of this work are summarized in Section~\ref{sec:Conclusion}, together with an outlook on future developments.

\section{Simulation environment and computational setup}
\label{sec:Simsetup}

Within this work, explicit finite element simulations are performed using \textsc{Radioss}. A user-defined material model is implemented that is automatically generated from a pretrained physics-augmented neural network (PANN). Before introducing the benchmark problem used to investigate different neural-network architectures and their influence on computational performance, the differences between implicit and explicit time integration schemes are briefly outlined, together with their implications for constitutive modeling.

Explicit solvers such as \textsc{Radioss} are widely used for highly dynamic problems, particularly crash, impact, and wave propagation simulations. For a general dynamic problem without damping, the governing equation of motion can be written as \cite{Wriggers.2008}
\begin{equation}
    \mathbf{M}\ddot{\mathbf{u}} + \mathbf{R\left(u\right)} = \mathbf{P},
    \label{eq:dynamic_equilibrium}
\end{equation}
where $\mathbf{M}$ denotes the mass matrix, $\mathbf{R}(\mathbf{u})$ the internal force vector, and $\mathbf{P}$ the external load vector.

In implicit time integration schemes, Eq.~(\ref{eq:dynamic_equilibrium}) is solved using an incremental-iterative procedure to determine the unknown displacement field at each time step. Typically, nonlinear algebraic systems are solved using Newton-type methods, which require repeated evaluations of the tangent stiffness matrix \cite{Wriggers.2008}. Consequently, the computational cost per time step is relatively high. On the other hand, implicit schemes can be formulated to be unconditionally stable, which makes them suitable for quasi-static or vibration-dominated problems.

In contrast, explicit formulations evaluate the solution in a forward manner, such that the state of the next time step depends only on quantities from previous time steps. Using the central difference scheme, the acceleration can be approximated as
\begin{equation}
    \ddot{\mathbf{u}} \approx \cfrac{\mathbf{u}_{n+1}-2\mathbf{u}_{n}+\mathbf{u}_{n-1}}{\left(\Delta t\right)^2},
    \label{eq:central_difference}
\end{equation}
where $\Delta t$ denotes the time step size. Substituting Eq.~(\ref{eq:central_difference}) into Eq.~(\ref{eq:dynamic_equilibrium}) yields
\begin{equation}
    \mathbf{M}\mathbf{u}_{n+1} = \left(\Delta t\right)^2\cdot
    \left(\mathbf{P}_n-\mathbf{R}\left(\mathbf{u}_n\right) \right)
    +\mathbf{M}\left(2\mathbf{u}_n-\mathbf{u}_{n-1}\right).
    \label{eq:explicit_update}
\end{equation}

This formulation can be evaluated efficiently without the need to compute the tangent stiffness matrix. For the central difference scheme, the admissible time-step size is governed by the largest numerical eigenfrequency $\omega_{\max}$ of the discretized system. Thus, the actual time step $\Delta t$ must be smaller than the cortical value $\Delta t_\mathrm{crit}$, which is defined as
\begin{equation}
\Delta t_{\mathrm{crit}} = \cfrac{2}{\omega_{\max}} .
\label{eq}
\end{equation}
Since an exact evaluation of $\omega_{\max}$ is generally impractical in large-scale finite element simulations, most explicit solvers employ a CFL-based estimate of the critical time step. This estimate relates the admissible time-step size to the time required for an elastic wave to pass through the smallest relevant finite element and is commonly written as \cite{Bathe.1996}
\begin{equation}
    \Delta t\approx\cfrac{h}{c_L}.
    \label{eq:critical_timestep}
\end{equation}
Here, $h$ denotes the characteristic element length and $c_L$ the longitudinal wave speed. For isotropic linear elasticity, the wave speed is given by
\begin{equation}
    c_L =  \sqrt{\cfrac{\lambda+2\mu}{\rho}},
    \label{eq:wavespeed_lame}
\end{equation}
with Lamé parameters $\lambda$ and $\mu$, and mass density $\rho$. Using the bulk modulus $\kappa = \lambda + \sfrac{2}{3}\,\mu$, Eq.~(\ref{eq:wavespeed_lame}) can alternatively be written as
\begin{equation}
    c_L =  \sqrt{\cfrac{3\kappa+4\mu}{3\rho}}.
    \label{eq:wavespeed_bulk}
\end{equation}

This relationship also forms the basis for \textit{mass scaling} \cite{Olovsson.2005, Tkachuk.2013}. In explicit simulations, the global stable time step is often determined by a small number of critical elements, such as very small, highly stiff, or strongly deformed elements. Instead of increasing the mass of the entire system, artificial mass is therefore typically added only locally to such critical elements in order to reduce their local wave speed and increase the corresponding critical time increment. This prevents isolated critical elements from disproportionately limiting the global time step.

Since explicit formulations do not require the tangent stiffness matrix, only the stress evaluation is considered in the following. Nevertheless, \textsc{Radioss} also supports implicit simulations, and the proposed implementation strategy can be extended to provide consistent tangent operators if required.

\subsection{Benchmark geometry and loading}
\label{sec:Setup_definition}

To generate a geometrically complex benchmark problem, an open-source dragon geometry obtained from \tovgu{MakerWorld}\footnote[1]{\href{https://makerworld.com/}{https://makerworld.com/}} was used. Since models designed for additive manufacturing frequently contain irregular surface triangulations and geometric artifacts, direct finite element meshing often leads to poor element quality and unstable simulations. To overcome these issues, the open-source software {\tovgu{MeshLab}}\footnote[2]{\href{https://www.meshlab.net/}{https://www.meshlab.net/}} was employed to preprocess the geometry.

Using the \textit{Isotropic Explicit Remeshing} tool available in MeshLab, clean and uniformly triangulated surface representations were generated for different target mesh resolutions. Besides improving the overall surface quality, the remeshing procedure also ensured a closed geometry, which is essential for reliable volume mesh generation. The resulting STL files are provided in the accompanying {\tovgu{Zenodo dataset}}\footnote[3]{\href{https://doi.org/10.5281/zenodo.20763660}{https://doi.org/10.5281/zenodo.20763660}}, allowing the benchmark geometries to be reproduced and used for further numerical studies. Three different surface discretizations were created and subsequently imported into \textsc{Simcenter HyperMesh}, the preprocessing environment provided by Siemens. There, tetrahedral volume meshes were generated. The resulting meshes are shown in Fig.~\ref{fig:Dragon_Mesh}.
\revI{
\begin{figure}[htb]
    \centering
    \begin{minipage}[t]{0.3\textwidth}
        \centering
        \textbf{Mesh 1}\\
        \vspace{1mm}
        {\raggedright
        \begin{tabular}{@{}l r@{}}
        Edge length:     & $1.6\,\mathrm{mm}$ \\
        No. of elements: & $85\,044$ \\
        No. of nodes:    & $17\,368$ \\
        DOFs:            & $52\,104$
        \end{tabular}
        \par}
        \vspace{1mm}
        \includegraphics[
            width=0.9\linewidth
        ]{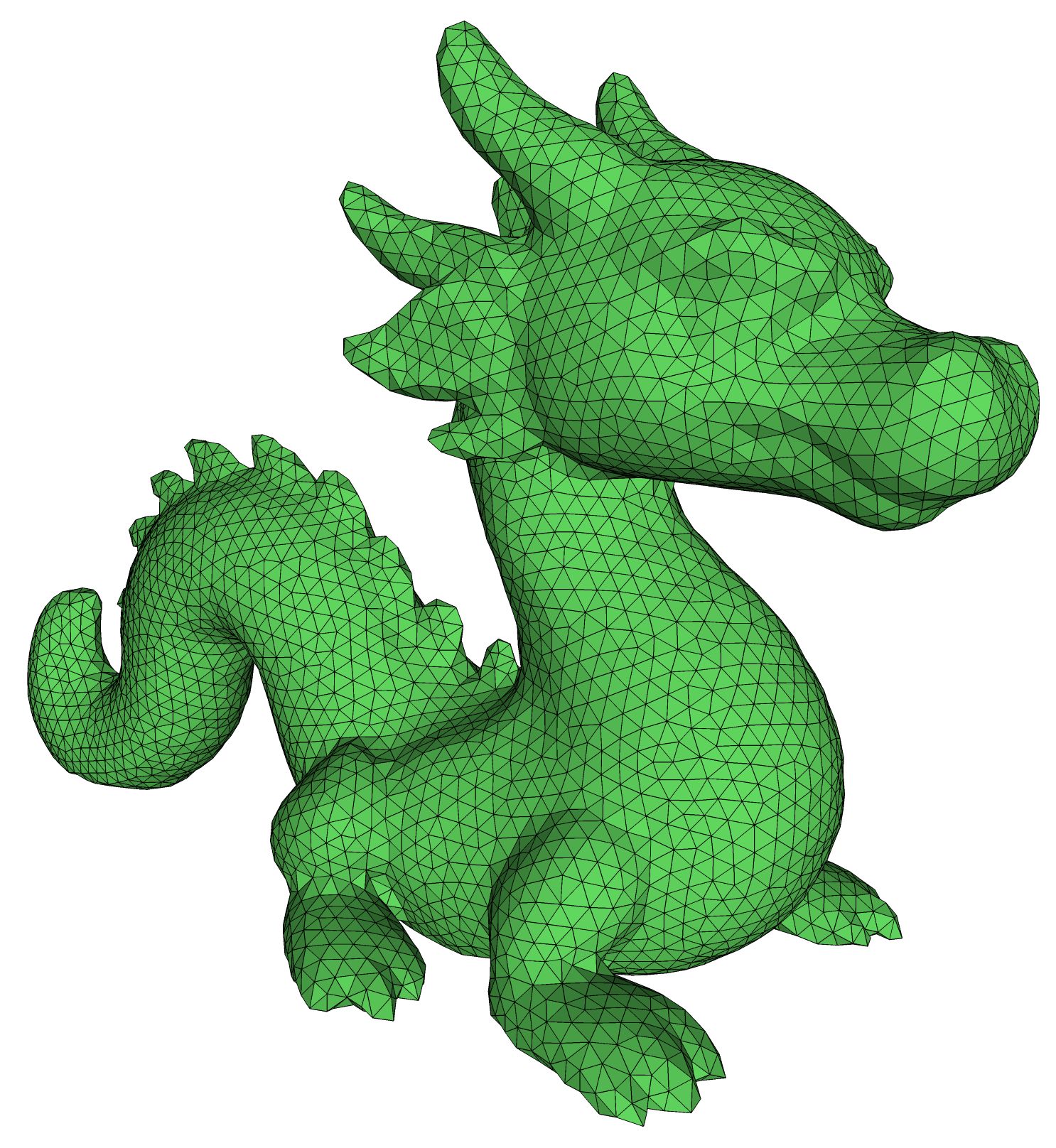}
    \end{minipage}
    \hfill
    \begin{minipage}[t]{0.3\textwidth}
        \centering
        \textbf{Mesh 2}\\
        \vspace{1mm}
        {\raggedright
        \begin{tabular}{@{}l r@{}}
        Edge length:     & $0.8\,\mathrm{mm}$ \\
        No. of elements: & $659\,143$ \\
        No. of nodes:    & $122\,305$ \\
        DOFs:            & $366\,915$
        \end{tabular}
        \par}
        \vspace{1mm}
        \includegraphics[
            width=0.9\linewidth
        ]{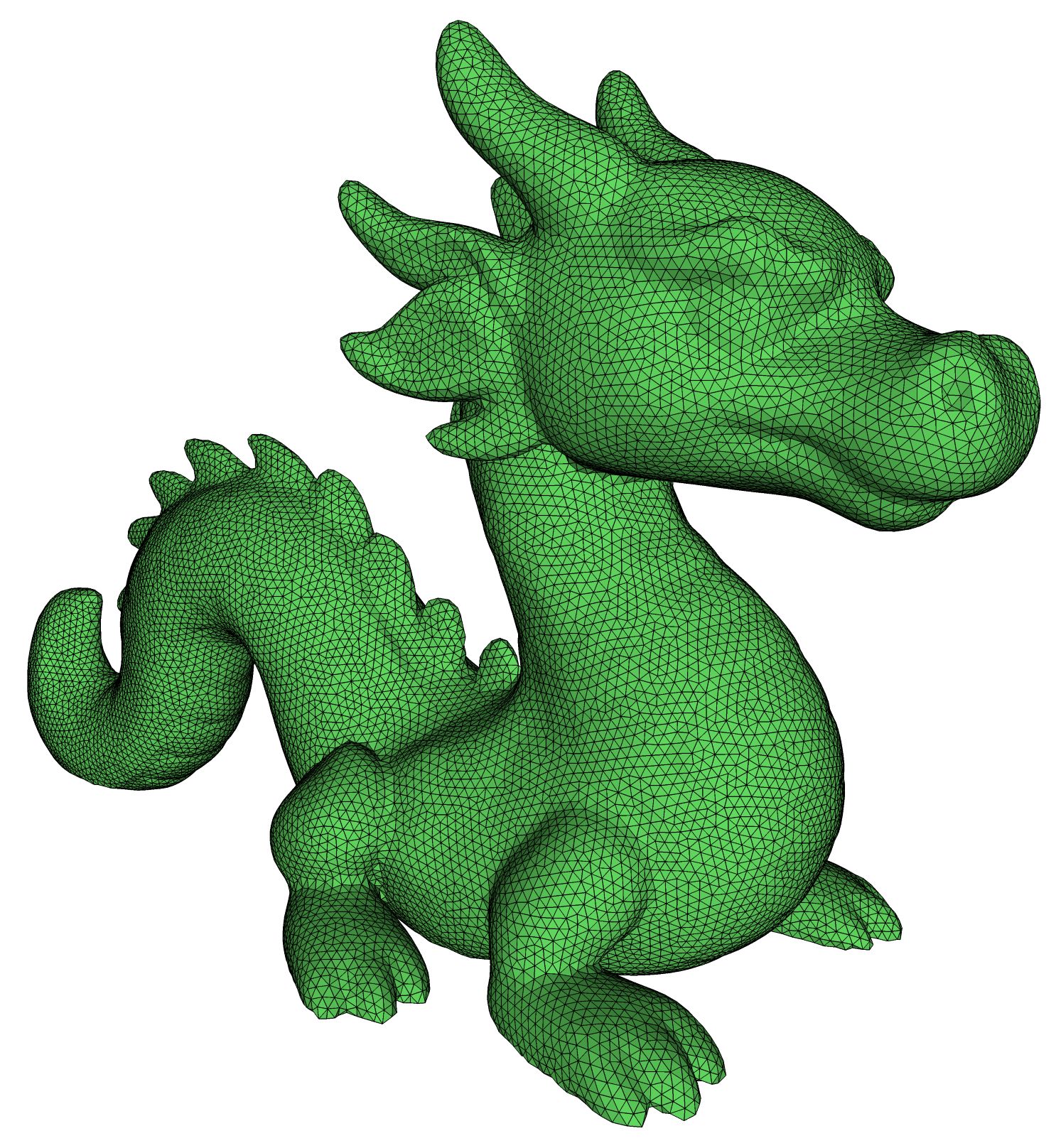}
    \end{minipage}
    \hfill
    \begin{minipage}[t]{0.3\textwidth}
        \centering
        \textbf{Mesh 3}\\
        \vspace{1mm}
        {\raggedright
        \begin{tabular}{@{}l r@{}}
        Edge length:     & $0.5\,\mathrm{mm}$ \\
        No. of elements: & $2\,669\,937$ \\
        No. of nodes:    & $475\,919$ \\
        DOFs:            & $1\,427\,757$
        \end{tabular}
        \par}
        \vspace{1mm}
        \includegraphics[
            width=0.9\linewidth
        ]{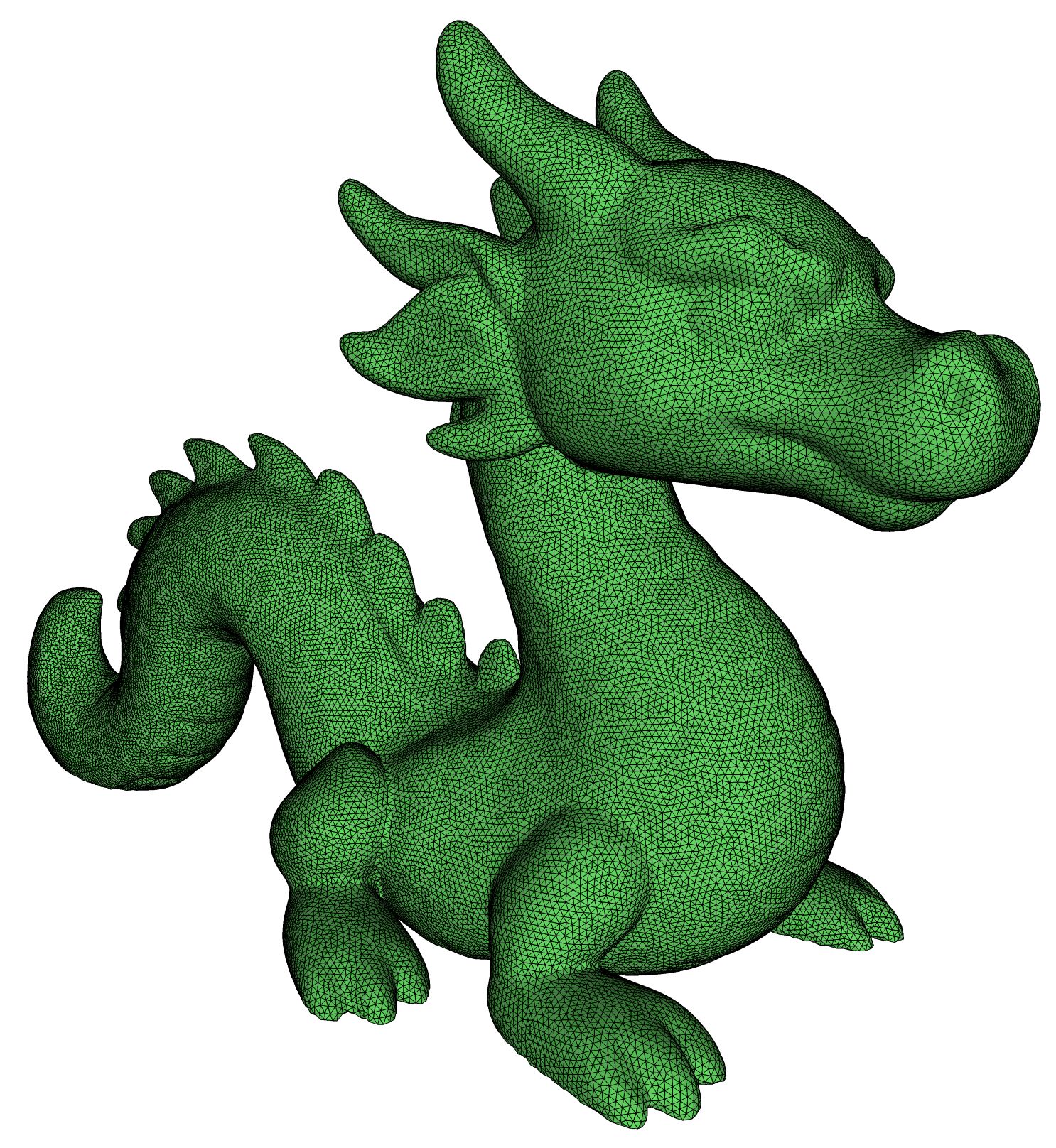}
    \end{minipage}
    \caption{Different mesh resolutions}
    \label{fig:Dragon_Mesh}
\end{figure}

The dragon model has an approximate length of}
\begin{equation}
l_{\mathrm{Dr}} \approx 10\,\mathrm{cm},
\end{equation}
and a volume of
\begin{equation}
V_{\mathrm{Dr}} \approx 54\,\mathrm{cm}^3.
\end{equation}

To evaluate the computational efficiency of the proposed constitutive formulations under highly dynamic loading conditions, an impact simulation is considered. A rigid spherical projectile with diameter
\begin{equation}
d_{\mathrm{Sp}} = 10\,\mathrm{mm}
\end{equation}
is launched toward the dragon model with an initial velocity of
\begin{equation}
v_{\mathrm{Sp}} = 101.8\,\mathrm{\cfrac{km}{h}}.
\end{equation}
The projectile is modeled as a rigid body with mass
\begin{equation}
m_{\mathrm{Sp}} = 4.1 \,\mathrm{g},
\end{equation}

which corresponds to the mass density of steel. The selected setup provides a representative explicit-dynamics benchmark, since the simulation contains large deformations, complex contact interactions, strongly nonlinear constitutive behavior, and a large number of evaluations of the material routine per time step. Consequently, the computational cost of the material model has a direct influence on the total simulation runtime, making this setup well suited for comparing classical constitutive formulations and PANN approaches.

\begin{figure}[htb]
    \centering
    \def\svgwidth{0.9\linewidth}
    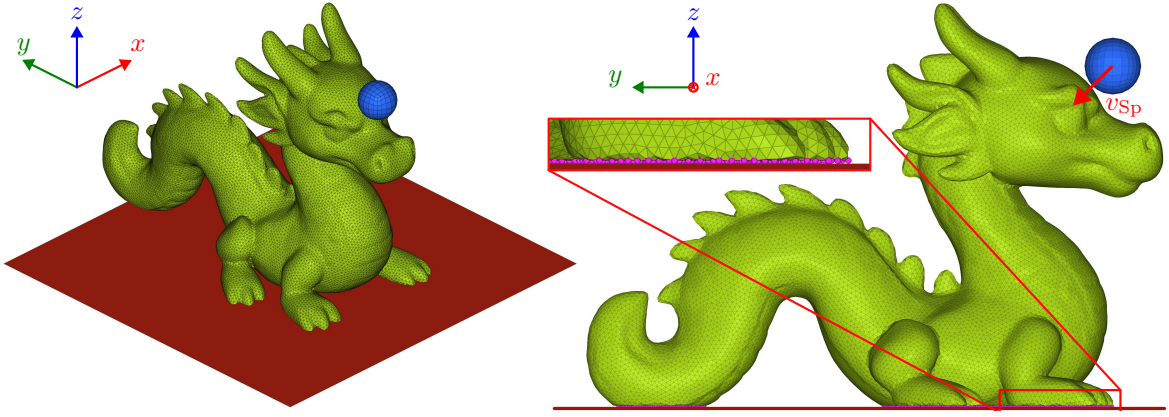
    \caption{Simulation setup and applied boundary conditions of the dragon model. The left image shows an isometric view of the discretized geometry on the rigid ground plane, while the right image shows a lateral view including the applied boundary conditions. The magenta nodes on the lower surface of the dragon are fixed in all translational directions.}
    \label{fig:Simulation_BCs}
\end{figure}

The dragon geometry was discretized using first-order tetrahedral elements. The nodes located on the lower surface of the dragon were fixed in all translational directions, while the remaining structure was allowed to deform freely. In addition, a rigid ground plane was defined below the dragon to prevent penetration of nodes through the lower boundary during the simulation. The simulation setup and the applied boundary conditions are shown in Fig.~\ref{fig:Simulation_BCs}.

\begin{figure}[htb]
    \centering
    \begin{minipage}[t]{0.3\textwidth}
        \centering
        \textbf{Mesh 1}\\
        \vspace{2mm}
        \setlength\tkwi{0.75\textwidth}
        \setlength\tkhe{0.8\tkwi}
\definecolor{mycolor1}{rgb}{0.00000,0.44700,0.74100}%
\begin{tikzpicture}

\begin{axis}[%
width=0.955\tkwi,
height=\tkhe,
at={(0\tkwi,0\tkhe)},
scale only axis,
xmin=0.26,
xmax=2.46,
xlabel style={font=\color{white!15!black}},
xlabel={Min. edge length $h_{\min}^{\mathrm{elem}}$ [mm]},
ymin=0,
ymax=25000,
ylabel style={font=\color{white!15!black}},
ylabel={Element count},
axis background/.style={fill=white},
xmajorgrids,
ymajorgrids
]
\addplot[ybar interval, fill=mycolor1, fill opacity=0.6, draw=black, area legend] table[row sep=crcr] {%
x	y\\
0.36	1\\
0.46	0\\
0.56	6\\
0.66	21\\
0.76	13\\
0.86	36\\
0.96	142\\
1.06	666\\
1.16	3066\\
1.26	9542\\
1.36	19428\\
1.46	23039\\
1.56	17384\\
1.66	7837\\
1.76	2744\\
1.86	860\\
1.96	210\\
2.06	40\\
2.16	7\\
2.26	2\\
2.36	2\\
};
\end{axis}
\end{tikzpicture}%
        \vspace{2mm}
        {\raggedright
        \begin{tabular}{@{}l r@{}}
        Mean $h_{\min}^{\mathrm{elem}}$: & $1.504\,\mathrm{mm}$\\[1mm]
        Min. $h_{\min}^{\mathrm{elem}}$: & $0.365\,\mathrm{mm}$\\[1mm]
        Max. $h_{\min}^{\mathrm{elem}}$: & $2.350\,\mathrm{mm}$ 
        \end{tabular}
        \par}
    \end{minipage}
    \hfill
    \begin{minipage}[t]{0.3\textwidth}
        \centering
        \textbf{Mesh 2}\\
        \vspace{2mm}
        \setlength\tkwi{0.75\textwidth}
        \setlength\tkhe{0.8\tkwi}
\definecolor{mycolor1}{rgb}{0.00000,0.44700,0.74100}%
\begin{tikzpicture}

\begin{axis}[%
width=0.955\tkwi,
height=\tkhe,
at={(0\tkwi,0\tkhe)},
scale only axis,
xmin=0.233,
xmax=1.267,
xlabel style={font=\color{white!15!black}},
xlabel={Min. edge length $h_{\min}^{\mathrm{elem}}$ [mm]},
ymin=0,
ymax=180000,
ylabel style={font=\color{white!15!black}},
ylabel={Element count},
axis background/.style={fill=white},
xmajorgrids,
ymajorgrids
]
\addplot[ybar interval, fill=mycolor1, fill opacity=0.6, draw=black, area legend] table[row sep=crcr] {%
x	y\\
0.28	3\\
0.327	0\\
0.374	12\\
0.421	24\\
0.468	169\\
0.515	1301\\
0.562	7437\\
0.609	32158\\
0.656	89537\\
0.703	166963\\
0.75	168111\\
0.797	113949\\
0.844	51390\\
0.891	18895\\
0.938	6640\\
0.985	1947\\
1.032	501\\
1.079	92\\
1.126	12\\
1.173	2\\
1.22	2\\
};
\end{axis}
\end{tikzpicture}%
        \vspace{2mm}
        {\raggedright
        \begin{tabular}{@{}l r@{}}
        Mean $h_{\min}^{\mathrm{elem}}$: & $0.762\,\mathrm{mm}$\\[1mm]
        Min. $h_{\min}^{\mathrm{elem}}$: & $0.289\,\mathrm{mm}$\\[1mm]
        Max. $h_{\min}^{\mathrm{elem}}$: & $1.201\,\mathrm{mm}$
        \end{tabular}
        \par}
    \end{minipage}
    \hfill
    \begin{minipage}[t]{0.3\textwidth}
        \centering
        \textbf{Mesh 3}\\
        \vspace{2mm}
        \setlength\tkwi{0.75\textwidth}
        \setlength\tkhe{0.8\tkwi}
\definecolor{mycolor1}{rgb}{0.00000,0.44700,0.74100}%
\begin{tikzpicture}

\begin{axis}[%
width=0.955\tkwi,
height=\tkhe,
at={(0\tkwi,0\tkhe)},
scale only axis,
xmin=0.118,
xmax=0.822,
xlabel style={font=\color{white!15!black}},
xlabel={Min. edge length $h_{\min}^{\mathrm{elem}}$ [mm]},
ymin=0,
ymax=800000,
ylabel style={font=\color{white!15!black}},
ylabel={Element count},
axis background/.style={fill=white},
xmajorgrids,
ymajorgrids
]
\addplot[ybar interval, fill=mycolor1, fill opacity=0.6, draw=black, area legend] table[row sep=crcr] {%
x	y\\
0.15	3\\
0.182	0\\
0.214	3\\
0.246	36\\
0.278	366\\
0.31	4779\\
0.342	34838\\
0.374	173388\\
0.406	477422\\
0.438	710017\\
0.47	620522\\
0.502	393591\\
0.534	168383\\
0.566	60867\\
0.598	19429\\
0.63	5116\\
0.662	1004\\
0.694	152\\
0.726	17\\
0.758	4\\
0.79	4\\
};
\end{axis}
\end{tikzpicture}%
        \vspace{2mm}
        {\raggedright
        \begin{tabular}{@{}l r@{}}
        Mean $h_{\min}^{\mathrm{elem}}$: & $0.470\,\mathrm{mm}$\\[1mm]
        Min. $h_{\min}^{\mathrm{elem}}$: & $0.171\,\mathrm{mm}$\\[1mm]
        Max. $h_{\min}^{\mathrm{elem}}$: & $0.786\,\mathrm{mm}$
        \end{tabular}
        \par}
    \end{minipage}
    \caption{Histograms of the element-wise minimum edge length $h_{\min}^{\mathrm{elem}}$ for the three investigated mesh discretizations. The listed values denote the mean, minimum, and maximum of $h_{\min}^{\mathrm{elem}}$ over all elements.}
    \label{fig:Hist_length}
\end{figure}
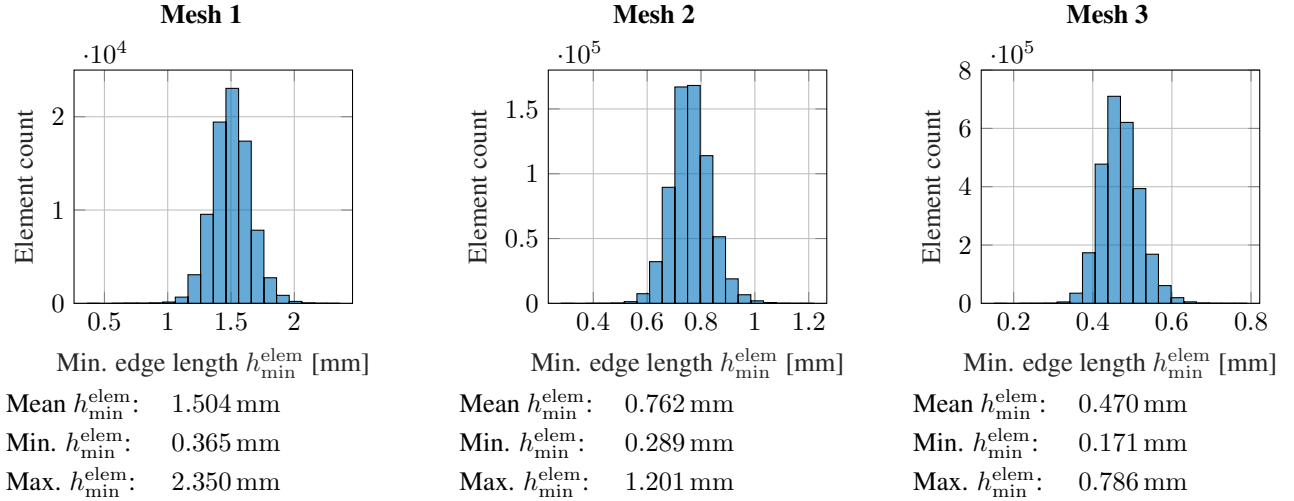

No mass scaling was applied. Despite the complex geometry of the dragon model, the generated meshes exhibited a sufficiently narrow distribution of the minimum element edge length $h_{\min}^{\mathrm{elem}}$, as shown in Fig.~\ref{fig:Hist_length}. This indicates that the element-size distribution did not introduce excessively small elements. Instead, the critical time increment $\Delta t_\mathrm{crit}$ was mainly governed by the contact interface. In \textsc{Radioss}, contact can reduce $\Delta t_\mathrm{crit}$ through the penalty stiffness and, for nonlinear contact formulations, through the kinematic condition that prevents nodes from crossing the contact gap within a single increment.

The contact treatment consisted of a self-contact definition for the deformable dragon structure and analytical rigid-wall definitions for the projectile and the ground plane. Self-contact was modeled using a penalty-based \texttt{/INTER/TYPE7} interface with a friction coefficient of $\mu=0.3$, a stiffness scale factor of $0.9$, and a minimum gap of $25\,\mathrm{\mu m}$. The rigid ground was represented by a rigid plane, \texttt{/RWALL/PLANE}, with a secondary-node search distance of $10\,\mathrm{mm}$ and a friction coefficient of $\mu=0.2$. The projectile was modeled as a spherical rigid wall \texttt{/RWALL/SPHER} with a mass of $4.1\,\mathrm{g}$ and a friction coefficient of $\mu=0.2$. Its initial velocity components were prescribed as $v_x=0$, $v_y=20\,\mathrm{m/s}$, and $v_z=-20\,\mathrm{m/s}$. The total simulated physical time was set to $t_\mathrm{end}=5\,\mathrm{ms}$. Time-history and three-dimensional simulation result files were written at intervals of $0.02\,\mathrm{ms}$. The requested output variables included the von Mises stress, nodal velocities, and the nodal added mass field, which was used to verify that no artificial mass was introduced during the simulations.

\subsection{Computational hardware}
\label{sec:CompHardware}

All finite element simulations were performed on the same workstation in order to ensure a consistent runtime comparison between the different material formulations. The hardware configuration used for all simulations is summarized in Table~\ref{tab:HardwareConf}.

\begin{table}[htb]
\centering
\caption{Hardware configuration used for the finite element simulations.}
\label{tab:HardwareConf}
\begin{tabular}{ll}
\hline
Component & Specification \\
\hline
CPU & Intel Core i9-13900K, 24 cores \\
Memory & 128 GB RAM \\
GPU & NVIDIA RTX A5500 \\
Operating system & Windows 11 Pro \\
Solver version & Siemens \textsc{Simcenter Radioss}, 2026 release \\
Parallel execution & 8 CPU threads, \texttt{-nt 8} \\
\hline
\end{tabular}
\end{table}

The reported finite element runtimes refer to the CPU-based \textsc{Radioss} simulations. The GPU was used only for neural-network training and visualization.
\section{Material model}
\label{sec:Matmodel}
While the projectile is modeled as a rigid body, the dragon structure is represented by a nearly incompressible rubber-like hyperelastic material. The objective of the present work is not the identification of a specific material for the benchmark geometry, but the integration and computational evaluation of PANN constitutive models in an explicit finite element framework. Therefore, the material behavior is restricted to isotropic hyperelasticity.

\begin{figure}[htb]
    \centering
    \def\svgwidth{0.6\linewidth}
    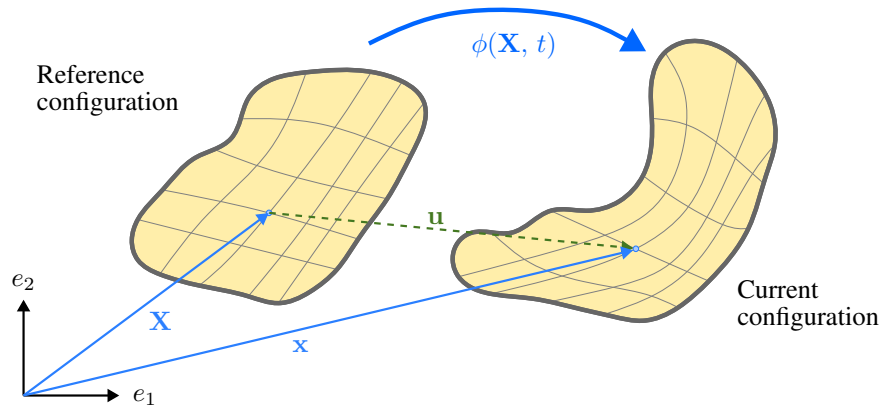
    \caption{Reference and current configurations of a deformable body. A material point with position vector $\mathbf{X}$ in the reference configuration is mapped by the motion $\boldsymbol{\phi}(\mathbf{X},t)$ to its spatial position $\mathbf{x}$ in the current configuration. The displacement vector is given by $\mathbf{u}=\mathbf{x}-\mathbf{X}$.}
    \label{fig:Configuration}
\end{figure}

Hyperelastic constitutive models are formulated in terms of a scalar-valued strain-energy density function $\psi$, which depends on the deformation state. Let $\mathbf{X}$ denote the position of a material point in the undeformed reference configuration and let \revI{ $\phi\left(\mathbf{X},t\right)$} describe its motion \cite{Belytschko.2014}. The reference and current configurations are then given by
\begin{equation}
    \mathbf{X}=\phi\left(\mathbf{X}, \, 0\right),
\end{equation}
and
\begin{equation}
    \mathbf{x}=\phi\left(\mathbf{X}, \, t\right).
\end{equation}
The relation between the reference configuration, the current configuration, and the motion $\boldsymbol{\phi}(\mathbf{X},t)$ is illustrated in Fig.~\ref{fig:Configuration}.

The deformation gradient is defined as
\begin{equation}
    \mathbf{F}=\cfrac{\partial\phi}{\partial\mathbf{X}}.
\end{equation}
from which strain and deformation measures such as the right Cauchy-Green deformation tensor can be obtained:
\begin{equation}
    \mathbf{C}=\mathbf{F}^\mathrm{T}\mathbf{F}.
\end{equation}
The principal stretches $\lambda_1$, $\lambda_2$, and $\lambda_3$ are the square roots of the eigenvalues of $\mathbf{C}$. They are frequently used in hyperelastic constitutive laws, for example in the incompressible Ogden model \cite{Ogden.1972}
\begin{equation}
    \psi_\mathrm{Og}=\sum_{i=1}^N\cfrac{\mu_i}{\alpha_i}\,\left(\lambda_1^{\alpha_i} + \lambda_2^{\alpha_i} + \lambda_3^{\alpha_i}-3\right).
\end{equation}
Alternatively, many hyperelastic models are formulated in terms of the invariants of $\mathbf{C}$,
\begin{align}
    I_1 &= \mathrm{tr}\left(\mathbf{C}\right)= \lambda_1^2+\lambda_2^2+\lambda_3^2,\\
    I_2 &= \mathrm{tr}\left(\mathrm{cof}\left(\mathbf{C}\right)\right)
    =\cfrac{1}{2}\,\left[\mathrm{tr}\,\left(\mathbf{C}\right)^2-\mathrm{tr}\,\left(\mathbf{C}^2\right)\right]=\lambda_1^2\lambda_2^2+\lambda_2^2\lambda_3^2+\lambda_3^2\lambda_1^2,\\
    I_3 &= \mathrm{det}\left(\mathbf{C}\right) = J^2 =\lambda_1^2\lambda_2^2\lambda_3^2.
\end{align}

For nearly incompressible formulations, it is advantageous to separate volumetric and isochoric deformation contributions by means of the Flory decomposition \cite{Flory.1961}. The deformation gradient is decomposed according to
\begin{equation}
\mathbf{F}= \mathbf{F}_\mathrm{iso}\mathbf{F}_\mathrm{vol},
\end{equation}
with the isochoric and volumetric parts defined as
\begin{equation}
\mathbf{F}_\mathrm{iso} = J^{-1/3}\mathbf{F},
\qquad\text{and}\qquad
\mathbf{F}_\mathrm{vol} = J^{1/3}\identMat.
\end{equation}
Here, $\identMat$ denotes the identity matrix. This decomposition leads to the isochoric invariants
\begin{equation}
\overline{I}_1=J^{-2/3}I_1,
\qquad
\overline{I}_2=J^{-4/3}I_2.
\end{equation}

In finite element analysis, incompressibility can be treated either by a mixed formulation, such as a \textit{u-p} formulation \cite{Wriggers.2008}, or by introducing a volumetric penalty term. In the mixed formulation, an additional pressure field $p$ is introduced as a Lagrange multiplier to enforce the constraint $J=1$ \cite{Bathe.1996}. In the present work, the penalty approach is adopted, since it is directly compatible with the explicit formulation used in \textsc{Radioss}. The strain-energy density function is decomposed into an isochoric and a volumetric contribution:
\begin{equation}
    \psi = \psi_\mathrm{iso}\left(\overline{I}_1,\,\overline{I}_2\right) + \psi_\mathrm{vol}\left(J\right).
\end{equation}
The isochoric part describes distortional deformation and includes models such as the Neo-Hooke, Mooney-Rivlin \cite{Mooney.1940, Rivlin.1948}, Arruda-Boyce \cite{Arruda.1993}, or Carroll \cite{Carroll.2011} models. The volumetric part depends only on the volume measure $J$. Common choices are the quadratic and logarithmic penalty terms \cite{Hartmann.2003}
\begin{align}
    \psi_\mathrm{vol}^\mathrm{quad} &= \cfrac{1}{2}\,\kappa_0 \left(J-1\right)^2,\\
    \psi_\mathrm{vol}^\mathrm{log} &= \cfrac{1}{2}\,\kappa_0 \left(\ln\left(J\right)\right)^2.
\end{align}
Here, $\kappa_0$ denotes the initial bulk modulus. Larger values of $\kappa_0$ enforce stronger resistance against volume changes and thus approximate the incompressible behavior more closely. However, a high ratio between bulk modulus and shear modulus may lead to volumetric locking, especially for standard low-order finite elements. Therefore, only nearly incompressible behavior is considered. For the purposes of this work, the decoupled material formulation introduced above is used as the basis for both the classical and neural-network-based constitutive models. To provide a representative calibration dataset containing different deformation modes, the classical experimental data by Treloar~\cite{Treloar.1944} are employed. This dataset includes uniaxial tension, equibiaxial tension, and pure shear loading conditions and is therefore well suited for comparing invariant-based hyperelastic material formulations. In the following, a classical phenomenological constitutive law is compared with a PANN formulation.

\subsection{Carroll material model}
\label{sec:Carroll}
As a classical reference model, the hyperelastic formulation proposed by Carroll \cite{Carroll.2011} is used. The model was developed to represent various experimental datasets for rubber-like materials, including the data by Treloar~\cite{Treloar.1944}. Its strain-energy density function is given by
\begin{equation}
	\psi = a\overline{I}_1 + b\overline{I}_1^4 + c\sqrt{\overline{I}_2} + \psi_\mathrm{vol},
\end{equation}
and the corresponding second Piola-Kirchhoff stress reads
\begin{equation}
	\mathbf{S} = 2\left(a+4b\overline{I}_1^3\right) \cfrac{\partial \overline{I}_1}{\partial \mathbf{C}} + c \cfrac{1}{\sqrt{\overline{I}_2}}\,\cfrac{\partial \overline{I}_2}{\partial \mathbf{C}} + 2\cfrac{\partial \psi_\mathrm{vol}}{\partial\mathbf{C}}.
\end{equation}
Here, the material parameters and documented data points from Steinmann et al.~\cite{Steinmann.2012} are used. The parameters are listed in Table~\ref{tab:Fit_Carroll}, and the fitted first Piola-Kirchhoff stress is shown in Fig.~\ref{fig:Fit_Carroll}.
\begin{table}[htb]
\begin{minipage}[c]{0.4\textwidth}
\centering
\caption{Material parameters for the Carroll model}
\vspace{0.5em}
\label{tab:Fit_Carroll}
\renewcommand{\arraystretch}{1.2}
\begin{tabular}{c|r l}
    $a$ & $0.15$ & $\mathrm{MPa}$\\
    \hline
    $b$ & $3.1\times10^{-7}$ & $\mathrm{MPa}$\\
    \hline
    $c$ & $0.095$ & $\mathrm{MPa}$
\end{tabular}
\end{minipage}
\hspace{2em}
\begin{minipage}[c]{0.5\textwidth}
\centering
\setlength\tkwi{0.85\textwidth}
\setlength\tkhe{0.6\tkwi}
\definecolor{mycolor1}{rgb}{0.00000,0.44706,0.74118}%
\definecolor{mycolor2}{rgb}{0.85098,0.32549,0.09804}%
\definecolor{mycolor3}{rgb}{0.46667,0.67451,0.18824}%
\begin{tikzpicture}

\begin{axis}[%
width=0.96\tkwi,
height=\tkhe,
at={(0\tkwi,0\tkhe)},
scale only axis,
xmin=1,
xmax=8,
xlabel style={font=\color{white!15!black}},
xlabel={$\lambda_1\;\left[-\right]$},
ymin=0,
ymax=7,
ylabel style={font=\color{white!15!black}},
ylabel={$P_{11}\;\left[\mathrm{MPa}\right]$},
axis background/.style={fill=white},
axis x line*=bottom,
axis y line*=left,
xmajorgrids,
ymajorgrids,
legend style={at={(0.14,-0.3)}, anchor=north west, legend columns=2,at={(0.03,0.97)}, anchor=north west, legend cell align=left, align=left, draw=white!15!black}
]
\addplot [color=mycolor1, dotted, mark size=0.7pt, mark=*, mark options={solid, mycolor1}]
  table[row sep=crcr]{%
1	0\\
1.01	0.0300000000000002\\
1.12	0.14\\
1.24	0.23\\
1.39	0.32\\
1.61	0.41\\
1.89	0.5\\
2.17	0.58\\
2.42	0.67\\
3.01	0.85\\
3.58	1.04\\
4.03	1.21\\
4.76	1.58\\
5.36	1.94\\
5.76	2.29\\
6.16	2.67\\
6.4	3.02\\
6.62	3.39\\
6.87	3.75\\
7.05	4.12\\
7.16	4.47\\
7.27	4.85\\
7.43	5.21\\
7.5	5.57\\
7.61	6.3\\
};
\addlegendentry{UT--Treloar}

\addplot [color=mycolor1]
  table[row sep=crcr]{%
1	0\\
1.01	0.0105261748323056\\
1.04666666666667	0.0471712158003426\\
1.08333333333333	0.0810667180658298\\
1.12	0.112575593163952\\
1.16	0.144586077074321\\
1.2	0.174456276343809\\
1.24	0.202468671330373\\
1.29	0.235231826807044\\
1.34	0.265851478958271\\
1.39	0.294646750782528\\
1.46333333333333	0.334119508274949\\
1.53666666666667	0.370891473360471\\
1.61	0.405476095427216\\
1.70333333333333	0.446951786995008\\
1.79666666666667	0.486159924997565\\
1.89	0.523576611979645\\
1.98333333333333	0.559565217480219\\
2.17	0.628325230121082\\
2.33666666666667	0.687130243315202\\
2.42	0.715884001021407\\
2.61666666666667	0.782560051153532\\
3.39	1.04066434378486\\
3.58	1.10591448257347\\
3.73	1.15875315720179\\
3.88	1.21314643477587\\
4.03	1.26949819963269\\
4.27333333333333	1.36628553052284\\
4.51666666666667	1.47180165545256\\
4.76	1.58893869322329\\
4.96	1.69639522455168\\
5.16	1.81644923577531\\
5.36	1.9518659748484\\
5.49333333333333	2.05224901036813\\
5.62666666666667	2.16195961142653\\
5.76	2.28218518607181\\
5.89333333333333	2.41423025010627\\
6.02666666666667	2.55952501253233\\
6.16	2.71963435813495\\
6.24	2.82351560504235\\
6.32	2.93373774450654\\
6.4	3.0507168756358\\
6.47333333333333	3.164256950805\\
6.54666666666667	3.28419583242113\\
6.62	3.41090343382047\\
6.70333333333333	3.56359713663252\\
6.78666666666667	3.72612657095604\\
6.87	3.8991170794747\\
6.93	4.03051517900855\\
6.99	4.16793161548324\\
7.05	4.31163090993103\\
7.08666666666667	4.40265699470747\\
7.12333333333333	4.49619570756257\\
7.16	4.59231259417572\\
7.19666666666667	4.69107456141006\\
7.23333333333333	4.79254989845205\\
7.27	4.89680829816798\\
7.32333333333333	5.05357151912786\\
7.37666666666667	5.21659796407585\\
7.45333333333333	5.46238648162773\\
7.5	5.61886911875211\\
7.53666666666667	5.74560025590441\\
7.57333333333333	5.87574773073298\\
7.61	6.00939530866806\\
};
\addlegendentry{UT--Model}

\addplot [color=mycolor2, dotted, mark size=0.7pt, mark=*, mark options={solid, mycolor2}]
  table[row sep=crcr]{%
1	0\\
1.04	0.0899999999999999\\
1.08	0.16\\
1.12	0.24\\
1.14	0.26\\
1.2	0.33\\
1.31	0.44\\
1.42	0.51\\
1.69	0.65\\
1.94	0.77\\
2.49	0.96\\
3.03	1.24\\
3.43	1.45\\
3.75	1.72\\
4.03	1.96\\
4.26	2.22\\
4.44	2.43\\
};
\addlegendentry{ET--Treloar}

\addplot [color=mycolor2]
  table[row sep=crcr]{%
1	0\\
1.02666666666667	0.0536566552795774\\
1.05333333333333	0.101730814115523\\
1.08	0.145092117033268\\
1.10666666666667	0.184448874055112\\
1.14	0.228888271750516\\
1.18	0.276422700258031\\
1.2	0.298188086163331\\
1.23666666666667	0.335159011650604\\
1.27333333333333	0.368878169103414\\
1.31	0.399875130783095\\
1.34666666666667	0.428581252379878\\
1.38333333333333	0.455351276605944\\
1.42	0.480479273955449\\
1.51	0.536651022264737\\
1.6	0.587050124178683\\
1.69	0.633495278673799\\
1.77333333333333	0.674050633150335\\
1.85666666666667	0.712954090564581\\
1.94	0.750696510767007\\
2.12333333333333	0.831304264169934\\
2.49	0.990240072574232\\
2.67	1.07026209573518\\
2.85	1.15347546971659\\
3.03	1.24148264382869\\
3.16333333333333	1.31088816235158\\
3.29666666666667	1.38491700221537\\
3.43	1.46466769178189\\
3.53666666666667	1.53344252984767\\
3.64333333333333	1.60744696734766\\
3.75	1.6875265655318\\
3.84333333333333	1.76332308747131\\
3.93666666666667	1.84520225377078\\
4.03	1.93394199333717\\
4.10666666666667	2.01256653909955\\
4.18333333333333	2.09691074559701\\
4.26	2.18752964636563\\
4.32	2.2632104479969\\
4.38	2.34340430195214\\
4.44	2.42843370803518\\
};
\addlegendentry{ET--Model}

\addplot [color=mycolor3, dotted, mark size=0.7pt, mark=*, mark options={solid, mycolor3}]
  table[row sep=crcr]{%
1	0\\
1.06	0.0700000000000003\\
1.14	0.16\\
1.21	0.24\\
1.32	0.33\\
1.46	0.42\\
1.87	0.59\\
2.4	0.76\\
2.98	0.93\\
3.48	1.11\\
3.96	1.28\\
4.36	1.46\\
4.69	1.62\\
4.96	1.79\\
};
\addlegendentry{PS--Treloar}

\addplot [color=mycolor3]
  table[row sep=crcr]{%
1	0\\
1.04	0.0535850704427716\\
1.06	0.0781893772162165\\
1.08666666666667	0.109008402206172\\
1.11333333333333	0.137790624351584\\
1.14	0.164759262550216\\
1.18666666666667	0.208160804135697\\
1.21	0.228276273490365\\
1.24666666666667	0.258035348993262\\
1.28333333333333	0.285797617041801\\
1.32	0.311822235087384\\
1.36666666666667	0.342774525689008\\
1.41333333333333	0.37164413678863\\
1.46	0.398742112725126\\
1.59666666666667	0.470149146315405\\
1.73333333333333	0.533012909855514\\
1.87	0.590176768375574\\
2.04666666666667	0.658473968137174\\
2.22333333333333	0.722645920911778\\
2.4	0.784200474238752\\
2.59333333333333	0.849726155113901\\
3.31333333333333	1.08993137747539\\
3.48	1.14696269640537\\
3.64	1.20300882352402\\
3.8	1.26074792091658\\
3.96	1.3206544619322\\
4.09333333333333	1.37262173817575\\
4.22666666666667	1.42682648202317\\
4.36	1.48366915340051\\
4.47	1.53286718089644\\
4.58	1.58444547658324\\
4.69	1.63871117652435\\
4.78	1.68534343326289\\
4.87	1.734199364131\\
4.96	1.78549397915\\
};
\addlegendentry{PS--Model}

\end{axis}
\end{tikzpicture}%

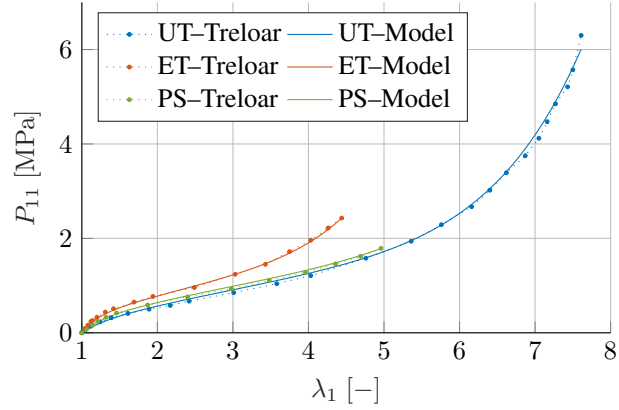
\captionof{figure}{Fitting result of the Carroll model}
\label{fig:Fit_Carroll}
\end{minipage}
\end{table}

\subsection{PANN model}
\label{sec:PANNmod}
The Carroll model serves as a reference for the proposed PANN formulation. In contrast to purely data-driven approaches, PANNs embed physical constraints directly into the architecture of the neural network. This improves robustness and helps to obtain constitutive responses that are consistent with basic principles of continuum mechanics.

For hyperelasticity, the neural network is used to represent the scalar strain-energy density function $\psi$ instead of predicting stresses directly. Stresses, and if required tangent stiffness matrices, can then be obtained by differentiating the energy with respect to the deformation measure. This construction ensures thermodynamic consistency of the stress response. In the present work, the deformation gradient $\mathbf{F}$ serves as the primary input to the constitutive routine. From this quantity, the right Cauchy-Green deformation tensor $\mathbf{C}$ and its invariants are computed and subsequently used as inputs to the neural network. This preserves objectivity, ensures symmetry of the stress tensor, and leads to an isotropic material response. Additional restrictions can be imposed through the network architecture. In particular, positive weights together with convex non-decreasing activation functions can be used to construct network representations that are both convex and monotonic with respect to the selected input invariants. A detailed discussion of such PANN architectures for hyperelasticity can be found in \cite{Linden.2023, Maurer.2024, Damma.2025b, Klein.2026}. In this paper, we follow the formulation proposed by Dammaß et al.~\cite{Damma.2025b}. Although the present study uses this specific PANN architecture, the implementation workflow can be easily adapted to alternative architectures.

For the construction of a polyconvex energy function, the choice of input invariants and the architecture of the neural network are both important \cite{Ebbing.2010}. Polyconvexity is a desirable mathematical property in hyperelasticity, since it is closely related to the existence of minimizers and contributes to stable and physically meaningful material responses under large deformations. By combining polyconvex input quantities with positive network weights and convex non-decreasing activation functions, the convexity properties of the invariants can be propagated through the network to the final strain-energy output. In this way, the resulting neural-network-based energy function can be constructed to preserve polyconvexity. While the first isochoric invariant $\overline{I}_1$ is polyconvex, the second isochoric invariant requires special treatment. As discussed by Klein et al.~\cite{Klein.2026}, a modified invariant
\begin{equation}
\overline{I}_2^\ast = \overline{I}_2^{3/2}
\end{equation}
can be used in this context. Therefore, the PANN formulation considered here uses $\overline{I}_1$ and $\overline{I}_2^\ast$ as input quantities.

The neural network used in the present work is an input-convex neural network (ICNN) \cite{Amos.2016}. In this class of networks, the weights are constrained to be positive and the activation functions are chosen to be convex and monotonically non-decreasing. This construction enables the strain-energy density to be represented as a convex function of the selected input invariants. The activation function most commonly used in this class of PANNs is the SoftPlus function,
\begin{equation}
    \mathcal{SP}\left(x\right) = \ln\left(1+\mathrm{e}^{x}\right).
\end{equation}
which is smooth, convex, and monotonically increasing. For a PANN with one hidden layer, the isochoric part of the strain-energy density can be written as
\begin{equation}
    \psi_\mathrm{iso} ^\mathrm{PANN} = \mathbf{w}^{(2)} \cdot \mathcal{SP}\left(\mathbf{w}^{(1)}\cdot\begin{bmatrix}\overline{I}_1 & \overline{I}_2^\ast\end{bmatrix}^\mathrm{T}+\boldsymbol{b}^{(1)}\right)
    \label{eq:Iso_PANN}
\end{equation}
Here, $\mathbf{w}^{(i)}$ denotes the weight matrix of layer $i$, and $\boldsymbol{b}^{(i)}$ denotes the corresponding bias vector. The first layer maps the selected invariant-based inputs to the hidden neurons, while the second layer combines the activated neuron responses into a scalar-valued strain-energy density. For the network architecture used in the present work, two invariant-based inputs, one hidden layer with five neurons, and one scalar energy output are employed. Consequently, the first weight matrix has dimensions $5\times2$, whereas the second one is of size $1\times5$.

\begin{figure}[htb]
    \centering
    \setlength\tkwi{0.425\textwidth}
    \setlength\tkhe{0.6\tkwi}
\definecolor{mycolor1}{rgb}{0.00000,0.44706,0.74118}%
\definecolor{mycolor2}{rgb}{0.85098,0.32549,0.09804}%
\definecolor{mycolor3}{rgb}{0.46667,0.67451,0.18824}%
\begin{tikzpicture}

\begin{axis}[%
width=0.96\tkwi,
height=\tkhe,
at={(0\tkwi,0\tkhe)},
scale only axis,
xmin=1,
xmax=8,
xlabel style={font=\color{white!15!black}},
xlabel={$\lambda_1\;\left[-\right]$},
ymin=0,
ymax=7,
ylabel style={font=\color{white!15!black}},
ylabel={$P_{11}\;\left[\mathrm{MPa}\right]$},
axis background/.style={fill=white},
axis x line*=bottom,
axis y line*=left,
xmajorgrids,
ymajorgrids,
legend style={at={(0.14,-0.3)}, anchor=north west, legend columns=2,at={(0.03,0.97)}, anchor=north west, legend cell align=left, align=left, draw=white!15!black}
]
\addplot [color=mycolor1, dotted, mark size=0.7pt, mark=*, mark options={solid, mycolor1}]
  table[row sep=crcr]{%
1	0\\
1.01	0.0300000000000002\\
1.12	0.14\\
1.24	0.23\\
1.39	0.32\\
1.61	0.41\\
1.89	0.5\\
2.17	0.58\\
2.42	0.67\\
3.01	0.85\\
3.58	1.04\\
4.03	1.21\\
4.76	1.58\\
5.36	1.94\\
5.76	2.29\\
6.16	2.67\\
6.4	3.02\\
6.62	3.39\\
6.87	3.75\\
7.05	4.12\\
7.16	4.47\\
7.27	4.85\\
7.43	5.21\\
7.5	5.57\\
7.61	6.3\\
};
\addlegendentry{UT--Treloar}

\addplot [color=mycolor1]
  table[row sep=crcr]{%
1	0\\
1.01	0.00994893709999989\\
1.12	0.108107874\\
1.24	0.197466067\\
1.39	0.292199247\\
1.61	0.410152048\\
1.89	0.539791011\\
2.17	0.656981887\\
2.42	0.755755113\\
3.01	0.978567269\\
3.58	1.19079514\\
4.03	1.36314419\\
4.76	1.67437475\\
5.36	2.00666676\\
5.76	2.30432246\\
6.16	2.68730683\\
6.4	2.97151631\\
6.62	3.28499837\\
6.87	3.73214257\\
7.05	4.13939393\\
7.16	4.43482736\\
7.27	4.77355439\\
7.43	5.35980201\\
7.5	5.65741148\\
7.61	6.18448477\\
};
\addlegendentry{UT--Model}

\addplot [color=mycolor2, dotted, mark size=0.7pt, mark=*, mark options={solid, mycolor2}]
  table[row sep=crcr]{%
1	0\\
1.04	0.0899999999999999\\
1.08	0.16\\
1.12	0.24\\
1.14	0.26\\
1.2	0.33\\
1.31	0.44\\
1.42	0.51\\
1.69	0.65\\
1.94	0.77\\
2.49	0.96\\
3.03	1.24\\
3.43	1.45\\
3.75	1.72\\
4.03	1.96\\
4.26	2.22\\
4.44	2.43\\
};
\addlegendentry{ET--Treloar}

\addplot [color=mycolor2]
  table[row sep=crcr]{%
1	0\\
1.04	0.0730527583000002\\
1.08	0.133828282\\
1.12	0.185190103\\
1.14	0.208028152\\
1.2	0.267650013\\
1.31	0.352783356\\
1.42	0.419203702\\
1.69	0.546733677\\
1.94	0.648397605\\
2.49	0.870459734\\
3.03	1.12508695\\
3.43	1.36841329\\
3.75	1.62610723\\
4.03	1.92375035\\
4.26	2.23401558\\
4.44	2.52715294\\
};
\addlegendentry{ET--Model}

\addplot [color=mycolor3, dotted, mark size=0.7pt, mark=*, mark options={solid, mycolor3}]
  table[row sep=crcr]{%
1	0\\
1.06	0.0700\\
1.14	0.16\\
1.21	0.24\\
1.32	0.33\\
1.46	0.42\\
1.87	0.59\\
2.4	0.76\\
2.98	0.93\\
3.48	1.11\\
3.96	1.28\\
4.36	1.46\\
4.69	1.62\\
4.96	1.79\\
};
\addlegendentry{PS--Treloar}

\addplot [color=mycolor3]
  table[row sep=crcr]{%
1	0\\
1.06	0.0738158845000001\\
1.14	0.155770348\\
1.21	0.216251788\\
1.32	0.296606566\\
1.46	0.381664508\\
1.87	0.576363022\\
2.4	0.783591455\\
3.48	1.17782379\\
3.96	1.35891629\\
4.36	1.52043406\\
4.69	1.66736553\\
4.96	1.80276122\\
};
\addlegendentry{PS--Model}

\end{axis}
\end{tikzpicture}%
    \caption{Fitting result of the PANN}
\label{fig:Fit_PANNcm}
\end{figure}
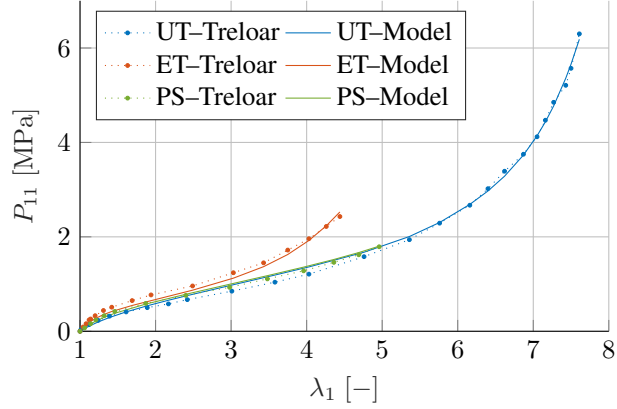

The PANN was trained several times using different numbers of layers and neurons. A sufficiently accurate fit was already obtained with one hidden layer containing five neurons. First-order Sobolev training with respect to the stress values was performed on the full dataset using the Sequential Least Squares Programming (SLSQP) algorithm. The resulting fit is shown in Fig.~\ref{fig:Fit_PANNcm}.

Both material descriptions reproduce the dataset well. In addition to the visual comparison in Figs.~\ref{fig:Fit_Carroll} and~\ref{fig:Fit_PANNcm}, the fit quality is evaluated using the mean squared error
\begin{equation}
    \mathrm{MSE} = \frac{1}{N}\sum_{j=1}^{N}\left(P_j^\mathrm{model}-P_j^\mathrm{data}\right)^2 .
\end{equation}
For the Carroll model and the PANN, the resulting errors are
\revI{
\begin{align}
    \mathrm{MSE}^\mathrm{Carroll} &= 4.22\times10^{-3}\,\mathrm{MPa}^2,\\
    \mathrm{MSE}^\mathrm{PANN} &= 5.07\times10^{-3}\,\mathrm{MPa}^2 .
\end{align}
}
Thus, the Carroll model performs slightly better for the present dataset. This can be attributed to the fact that the PANN formulation cannot reproduce the non-polyconvex term $\sqrt{\overline{I}_2}$ used in the Carroll model. In general, enforcing polyconvexity through a convex and monotonic formulation ensures stability by design, but it may restrict the permissible range of behavior more than necessary \cite{Klein.2026}. Nevertheless, both models provide a very good fit. In subsequent simulations, however, small differences in the material response are expected.
\section{Implementation}
\label{sec:Implementation}
For the implementation of the neural-network-based constitutive model, a user material subroutine is used in \textsc{Radioss}. The routine is written in Fortran and evaluates the Cauchy stress $\boldsymbol{\sigma}$ from the deformation gradient $\mathbf{F}$. To avoid a manual and error-prone transfer of trained network parameters into the finite element code, a Python script is used to automatically generate the corresponding material routine. The generated code contains the network architecture, the trained weights and biases, the activation functions, the invariant calculation, and the final stress evaluation. The script is provided in a public {\tovgu{GitHub repository}}\footnote[1]{\href{https://github.com/OVGU-CoMe/PANN_Radioss}{https://github.com/OVGU-CoMe/PANN\_Radioss}} and can be adapted to different network architectures.

In contrast to the training procedure and to the stress evaluation used for the fitted material curves, the standalone Fortran material routine cannot rely on automatic differentiation provided by PyTorch or TensorFlow. Wrapper-based approaches are possible in principle, but introduce additional dependencies and were found to be slower in our tests than a fully hard-coded stress evaluation. Therefore, the derivatives required for the stress calculation are explicitly generated by the Python script. This also enables a consistent runtime comparison between the classical Carroll model and the neural-network-based material model, since only the constitutive evaluation inside \textsc{Radioss} is considered.

After training, the Python script exports the trained PANN into a Fortran module. In the following, the stress formulation is derived for a PANN with one hidden layer. Symbols carrying a hat, e.g., $\hat{\mathbf{S}}$, denote quantities written in Voigt notation using the component ordering
\begin{equation}
    \hat{\mathbf{S}} = \begin{bmatrix}
        S_{11} & S_{22} & S_{33} & S_{12} & S_{23} & S_{13}  
    \end{bmatrix}^\mathrm{T}.
\end{equation}

The strain-energy density of the PANN from Eq.~\eqref{eq:Iso_PANN}, extended by the volumetric penalty term, is written as
\revI{
\begin{equation}
\begin{split}
    \psi^\mathrm{PANN} &= \psi_\mathrm{iso} ^\mathrm{PANN}\left(\overline{I}_1,\,\overline{I}_2^\ast\right)  + 
    \psi_\mathrm{vol}\left(J\right)\\
    &= \mathbf{w}^{(2)} \cdot \mathcal{SP}\left(\mathbf{w}^{(1)}\cdot\begin{bmatrix}\overline{I}_1 & \overline{I}_2^\ast\end{bmatrix}^\mathrm{T}+\boldsymbol{b}^{(1)}\right) + \cfrac{1}{2}\,\kappa_0\left(J-1\right)^2.
\end{split}
\label{eq:Energy_PANN}
\end{equation}
}
For this network structure, the second Piola-Kirchhoff stress follows from the first derivative of the strain-energy density with respect to the right Cauchy-Green deformation tensor. Applying the chain rule gives
\begin{equation}
\begin{split}
    \hat{\mathbf{S}}^\mathrm{PANN} & = 2\cfrac{\partial\psi_\mathrm{iso}^\mathrm{PANN}}{\partial\hat{\mathbf{C}}}+2\cfrac{\partial\psi_\mathrm{vol}}{\partial\hat{\mathbf{C}}}\\
    & =
    2\cdot \mathbf{w}^{(2)} \cdot \cfrac{\partial\, \mathcal{SP}\left(\mathbf{w}^{(1)}\cdot\begin{bmatrix}\overline{I}_1 & \overline{I}_2^\ast\end{bmatrix}^\mathrm{T}+\boldsymbol{b}^{(1)}\right)}{\partial\hat{\mathbf{C}}}
    +2\cfrac{\partial\psi_\mathrm{vol}}{\partial\hat{\mathbf{C}}}\\
    & =
    2\cdot \mathbf{w}^{(2)} \cdot \mathrm{diag}\left(\mathcal{SP}'\left(\dots\right)\right) \cdot \mathbf{w}^{(1)} \cdot 
    \cfrac{\partial \begin{bmatrix}\overline{I}_1 & \overline{I}_2^\ast\end{bmatrix}^\mathrm{T}}{\partial\hat{\mathbf{C}}}
    +2\cfrac{\partial\psi_\mathrm{vol}}{\partial\hat{\mathbf{C}}}
    \\
    & = 
    2\cdot \mathbf{w}^{(2)} \cdot \mathrm{diag}\left(\mathcal{SP}'\left(\dots\right)\right) \cdot \mathbf{w}^{(1)}\cdot
    \begin{bmatrix}
        \cfrac{\partial \overline{I}_1}{\partial C_{11}} & \cfrac{\partial \overline{I}_1}{\partial C_{22}} & \cfrac{\partial \overline{I}_1}{\partial C_{33}} & \cfrac{\partial \overline{I}_1}{\partial C_{12}} & \cfrac{\partial \overline{I}_1}{\partial C_{23}} & \cfrac{\partial \overline{I}_1}{\partial C_{13}} \\[12pt]
        \cfrac{\partial \overline{I}_2^\ast}{\partial C_{11}} & \cfrac{\partial \overline{I}_2^\ast}{\partial C_{22}} & \cfrac{\partial \overline{I}_2^\ast}{\partial C_{33}} & \cfrac{\partial \overline{I}_2^\ast}{\partial C_{12}} & \cfrac{\partial \overline{I}_2^\ast}{\partial C_{23}} & \cfrac{\partial \overline{I}_2^\ast}{\partial C_{13}}
    \end{bmatrix}
    +2\cfrac{\partial\psi_\mathrm{vol}}{\partial\hat{\mathbf{C}}}.
    \label{eq:PK2_PANN}
\end{split}
\end{equation}
For one hidden layer, this expression already contains the complete network derivative. For additional hidden layers, one additional weight matrix and one diagonal matrix containing the activation derivatives are inserted for each layer.

The derivatives of the isochoric invariants with respect to $\mathbf{C}$ are
\begin{align}
    \cfrac{\partial \overline{I}_1}{\partial \mathbf{C}} &= J^{-2/3} \left(\identMat - \cfrac{1}{3}\,I_1 \mathbf{C}^{-1}\right),\\
    \cfrac{\partial \overline{I}_2}{\partial \mathbf{C}} &= J^{-4/3}\left(I_1\identMat - \mathbf{C} - \cfrac{2}{3}\,I_2 \mathbf{C}^{-1}\right).
\end{align}
If the modified invariant $\overline{I}_2^\ast$ is used, its derivative follows from an additional chain-rule factor:
\begin{equation}
    \cfrac{\partial \overline{I}_2^\ast}{\partial \mathbf{C}}
    =\cfrac{\partial \overline{I}_2^{3/2}}{\partial \mathbf{C}}
    =
    \cfrac{3}{2}\,
    \overline{I}_2^{1/2}
    \cfrac{\partial \overline{I}_2}{\partial \mathbf{C}}.
\end{equation}

For the volumetric part, the quadratic penalty formulation is used,
\begin{equation}
    \psi_\mathrm{vol}
    =
    \cfrac{1}{2}\,\kappa_0\left(J-1\right)^2 .
\end{equation}
This choice is consistent with the nearly incompressible formulation introduced in Section~\ref{sec:Matmodel}. Since the simulations are performed close to the incompressible regime, the specific choice of the volumetric penalty function is not expected to dominate the comparison between the material models. Its contribution to the second Piola-Kirchhoff stress is
\begin{equation}
\begin{split}
    2\cfrac{\partial\psi_\mathrm{vol}}{\partial\mathbf{C}} &=
    2\cfrac{\partial\left(\cfrac{1}{2}\kappa_0\left(J-1\right)^2\right)}{\partial\mathbf{C}}\\
    & = 2\kappa_0 \left(J-1\right) \cfrac{\partial J}{\partial\mathbf{C}}.
\end{split}   
\end{equation}
The derivative of $J$ with respect to the right Cauchy-Green deformation tensor is
\begin{equation}
    \cfrac{\partial J}{\partial \mathbf{C}} = \cfrac{1}{2}\, J \mathbf{C}^{-1}.
\end{equation}

In a direct implementation based on $\mathbf{C}$, the second Piola-Kirchhoff stress would finally have to be transformed into the Cauchy stress by
\begin{equation}
    \boldsymbol{\sigma} = \cfrac{1}{J}\,\mathbf{F}\mathbf{S}\mathbf{F}^\mathrm{T}.
\end{equation}
Since this operation is performed at every integration point and every time step, the generated user material evaluates the stress directly using the left Cauchy-Green deformation tensor
\begin{equation}
    \mathbf{b} = \mathbf{F}\mathbf{F}^\mathrm{T}.
\end{equation}
The invariants of $\mathbf{b}$ are identical to those of $\mathbf{C}$, so the definition of the PANN remains unchanged. The Kirchhoff stress can then be evaluated directly as
\begin{equation}
    \boldsymbol{\tau} = J\boldsymbol{\sigma} = 2\cfrac{\partial\psi}{\partial\mathbf{b}}\,\mathbf{b}.
\end{equation}
The derivatives required for this formulation can be written in a form where they are already multiplied by $\mathbf{b}$:
\begin{align}
    \cfrac{\partial \overline{I}_1}{\partial \mathbf{b}}\,\mathbf{b} &= J^{-2/3} \left(\mathbf{b} - \cfrac{1}{3}\,I_1 \identMat\right),\label{eq:Deriv_bmat1}\\
    \cfrac{\partial \overline{I}_2}{\partial \mathbf{b}}\,\mathbf{b} &= J^{-4/3}\left(I_1\mathbf{b} - \mathbf{b}^2 - \cfrac{2}{3}\,I_2 \identMat\right),\label{eq:Deriv_bmat2}\\
    \cfrac{\partial J}{\partial \mathbf{b}}\,\mathbf{b} &= \cfrac{1}{2}\, J \identMat.
    \label{eq:Deriv_bmat3}
\end{align}
With these relations, the complete stress evaluation can be written in a compact form and translated into the generated Fortran material module.

\subsection{Parameter definition and activation function}
\label{sec:Implem_Activation}

Since the PANN architecture may vary between trained models, the generated routine starts with the definition of the network size, the trained weights and biases, and the material parameters. The initial shear modulus of the network is also exported. This value is required to choose a suitable bulk modulus and to provide the wave-speed information used by \textsc{Radioss}, see Eq.~\eqref{eq:wavespeed_bulk}. The initial shear modulus depends on the isochoric part of the strain-energy density and is calculated as \cite{Holzapfel.2010}
\begin{equation}
    \mu_0 = 2\left(\cfrac{\partial \psi}{\partial I_1}+\cfrac{\partial \psi}{\partial I_2}\right)\Bigg|_{I_1=3,\, I_2=3,\, J=1}.
\end{equation}
For a network using $\overline{I}_2^*$ as input, the derivative with respect to $\overline{I}_2$ includes the corresponding chain-rule factor
\begin{equation}
\cfrac{\partial \psi}{\partial \overline{I}_2}
=
\cfrac{\partial \psi}{\partial \overline{I}_2^\ast}
\cfrac{\partial \overline{I}_2^\ast}{\partial \overline{I}_2}
=
\cfrac{3}{2}\sqrt{\overline{I}_2}\,
\cfrac{\partial \psi}{\partial \overline{I}_2^\ast}.
\end{equation}
Evaluated in the reference configuration, this gives
\begin{equation}
\mu_0
=
2\left(
\cfrac{\partial \psi}{\partial \overline{I}_1}
+
\cfrac{3}{2}\sqrt{3}\,
\cfrac{\partial \psi}{\partial \overline{I}_2^\ast}
\right)\Bigg|_{\overline{I}_1=3,\, \overline{I}_2=3,\, J=1}.
\end{equation}
To represent nearly incompressible behavior, the bulk modulus is chosen relative to the initial shear modulus. In the present simulations, we use approximately $\kappa_0 \approx 90\,\mu_0$. For the pretrained PANN considered here, this gives
\begin{align}
    \mu_0^\mathrm{PANN} &= 0.334\,\mathrm{MPa},\\
    \kappa_0^\mathrm{PANN} &= 30\,\mathrm{MPa}.
\end{align}
To generate the material module, the Python script loads the main quantities listed in Table~\ref{table:Load_variables}.
\begin{table}[htbp]
    \centering
    \caption{Overview of the main variables used for code generation}
    \label{table:Load_variables}
    \begin{tabular}{lll}
    \hline
    \textbf{Variable} & \textbf{Function} & \textbf{Structure} \\
    \hline
    \texttt{layers} &
    Stores the neural network layers and their parameters &
    List of layer objects or tuples \\
    \texttt{sizes} &
    Defines the number of neurons in each layer &
    Integer list, e.g. \texttt{[2, 4, 1]} \\
    \texttt{L} &
    Number of network layers &
    Integer scalar \\
    \texttt{k0} &
    Initial bulk modulus used for the volumetric energy contribution &
    Floating-point scalar \\
    \texttt{init\_shear} &
    Initial shear modulus used for material characterization &
    Floating-point scalar \\
    \hline
    \end{tabular}
\end{table}

The generated material routine begins with the module header. This header defines the numerical precision, the network architecture, the weight and bias arrays, and the material parameters required for the constitutive evaluation:
\begin{lstlisting}[
    style=fortranstyle,
    caption={Header of the generated Fortran material module.},
    label={lst:pann_header}
]
module pann_model_001
  implicit none
  integer, parameter :: dp = selected_real_kind(15, 307)   ! double precision kind

  logical, save :: params_loaded = .false.   ! flag indicating whether parameters have been loaded

  integer, parameter :: NL = 2   ! number of network layers
  integer, parameter :: N0 = 2   ! number of input neurons / invariants
  integer, parameter :: N1 = 5   ! number of hidden neurons
  integer, parameter :: N2 = 1   ! number of output neurons / energy output

  real(dp) :: W1_MAT(5,2)  ! weight matrix of Layer 1
  real(dp) :: b1_VEC(5)   ! bias vector of Layer 1

  real(dp) :: W2_MAT(1,5)  ! weight matrix of Layer 2
  ! Layer 2: no bias stored

  real(dp) :: K   ! bulk modulus
  real(dp) :: G   ! shear modulus

contains
\end{lstlisting}

The activation function is defined next. Two routines are generated: one routine evaluates both the SoftPlus function and its derivative, while the second routine evaluates only the derivative. The function value is required for networks with more than one hidden layer, whereas the derivative is required in all cases for the stress calculation.

For numerical robustness, the SoftPlus function is evaluated in a stabilized form. Instead of computing $\log(1+\exp(x))$ directly, which may overflow for large positive values of $x$, the SoftPlus function can be written as the two-term log-sum-exp expression
\begin{align}
    \mathcal{SP}(x) &=  \log(1+\exp(x))\\
                    &=  \log\left(\exp(0)+\exp(x)\right).
\end{align}
Following the shifted log-sum-exp formulation \cite{Blanchard.2021}, this expression can be rewritten as
\begin{equation}
    \mathcal{SP}(x)
    =
    m+\log\left(\exp(-m)+\exp(x-m)\right),
\end{equation}
where $m$ is a real-valued shift. Choosing $m=\max(x,0)$ prevents overflow in the exponential terms and gives
\begin{equation}
\mathcal{SP}(x)
=
\begin{cases}
\ln\left(1+\exp(x)\right), & x \leq 0, \\[0.3em]
x + \ln\left(1+\exp(-x)\right), & x > 0.
\end{cases}
\end{equation}
Equivalently, this stabilized expression can be written in compact form as
\begin{equation}
\mathcal{SP}(x)
=
\max(x,0) + \ln\left(1+\exp\left(-\left|x\right|\right)\right).
\end{equation}
This formulation only evaluates exponentials of non-positive arguments and is therefore robust for both large positive and large negative inputs. The same stabilization is used for the derivative to avoid overflow in the exponential terms. The Fortran implementation contains two variants of the activation routine. The first routine evaluates both the SoftPlus function and its derivative. The second routine evaluates only the derivative of the activation function and is used in parts of the material routine where the function value itself is not required. This avoids unnecessary function evaluations and reduces the computational cost. The corresponding implementations are shown in Listings~\ref{lst:pann_activation_SP} and~\ref{lst:pann_activation_SPD}.
\revI{
\begin{lstlisting}[
    style=fortranstyle,
    caption={SoftPlus activation function and derivative of the generated Fortran material module.},
    label={lst:pann_activation_SP},
    firstnumber=23,
]
  pure elemental subroutine softplus(x, y, dy)
    real(dp), intent(in)  :: x
    real(dp), intent(out) :: y, dy
    real(dp) :: ax, t
      ax = abs(x)
      t  = exp(-ax)
      y  = max(x, 0.0_dp) + log(1.0_dp + t)
      dy = merge( 1.0_dp/(1.0_dp+t), t/(1.0_dp+t), x >= 0.0_dp )
  end subroutine softplus
\end{lstlisting}
\begin{lstlisting}[
    style=fortranstyle,
    caption={Derivative-only SoftPlus activation routine of the generated Fortran material module.},
    label={lst:pann_activation_SPD},
    firstnumber=42,
]
  pure elemental subroutine softplusD(x, dy)
    real(dp), intent(in)  :: x
    real(dp), intent(out) :: dy
    real(dp) :: ax, t
      ax = abs(x)
      t  = exp(-ax)
      dy = merge( 1.0_dp/(1.0_dp+t), t/(1.0_dp+t), x >= 0.0_dp )
  end subroutine softplusD
\end{lstlisting}
}
After defining the activation routines, the trained PANN parameters are assigned in the subroutine \texttt{load\_nn\_params()}. This routine defines the material parameter values, normalization constants, and neural-network weights used by the generated material model, as shown in Listing~\ref{lst:pann_matparams}.
\begin{lstlisting}[
    style=fortranstyle,
    caption={Parameter-loading routine of the generated Fortran material module.},
    label={lst:pann_matparams},
    firstnumber=62,
]
  subroutine load_nn_params()
    implicit none

    K = 3.00000000000000000e+01_dp
    G = 3.34615804329204680e-01_dp

    ! -----------------------------
    ! Layer 1: W1_MAT(5,2)
    ! -----------------------------
    W1_MAT(1,1) = 9.59297979487586279e-03_dp
    W1_MAT(1,2) = 4.24461749664615795e-07_dp
    W1_MAT(2,1) = 2.85353781885561819e-01_dp
    W1_MAT(2,2) = 6.81215943898273779e-01_dp
    W1_MAT(3,1) = 2.04024153429259780e-01_dp
    W1_MAT(3,2) = 3.27018085970901829e-08_dp
    W1_MAT(4,1) = 9.18260524707273512e-02_dp
    W1_MAT(4,2) = 4.36300720152440421e-12_dp
    W1_MAT(5,1) = 2.93601395665652409e-01_dp
    W1_MAT(5,2) = 2.58865649108586606e-05_dp

    ! Layer 1: b1_VEC(5)
    b1_VEC(1) = -5.97115873111896178e-01_dp
    b1_VEC(2) = -5.58964695851211446e-01_dp
    b1_VEC(3) = -6.81407455398122686e+00_dp
    b1_VEC(4) = -7.88836723395510742e+00_dp
    b1_VEC(5) = 1.84126502958735685e+01_dp

    ! -----------------------------
    ! Layer 2: W2_MAT(1,5)
    ! -----------------------------
    W2_MAT(1,1) = 1.27447829295857674e-02_dp
    W2_MAT(1,2) = 1.30974123113002333e-04_dp
    W2_MAT(1,3) = 1.04629411280879783e-01_dp
    W2_MAT(1,4) = 3.30124886772518025e+01_dp
    W2_MAT(1,5) = 5.63979314360676831e-01_dp

  end subroutine load_nn_params
\end{lstlisting}

To avoid repeated assignments during the simulation, the separate routine \texttt{ensure\_nn\_params\_loaded()} checks whether the parameters have already been initialized. Thus, the parameter assignment is performed only once. The generated parameter block is shown in Listing~\ref{lst:pann_ensure_loaded}.
\begin{lstlisting}[
    style=fortranstyle,
    caption={Material parameters of the generated PANN model.},
    label={lst:pann_ensure_loaded},
    firstnumber=100,
]
  subroutine ensure_nn_params_loaded()
    if (.not. params_loaded) then
      call load_nn_params()
      params_loaded = .true.
    end if
  end subroutine
\end{lstlisting}

At this point, the network architecture, the activation functions, and the material parameters are available. The following section describes how the user material evaluates the constitutive response: first, the required kinematic quantities are computed from the deformation gradient; subsequently, the trained PANN model is used to obtain the network derivatives and the resulting Cauchy stress.

\subsection{Constitutive quantities evaluated in the user material}

The stress evaluation is formulated in terms of the left Cauchy-Green deformation tensor. The deformation gradient passed to the material routine by \textsc{Radioss} is stored in the vector format
\begin{equation}
    \mathbf{F} = \begin{bmatrix}
        F_{11} & F_{22} & F_{33} & F_{12} & F_{23} & F_{13} & F_{21} & F_{32} & F_{31}
    \end{bmatrix}.
\end{equation}
From this input, the tensor $\mathbf{b}$ and its square $\mathbf{b}^2$ are computed in Voigt notation. Subsequently, the invariants and the invariant derivatives multiplied by $\mathbf{b}$ are evaluated according to Eqs.~\eqref{eq:Deriv_bmat1}--\eqref{eq:Deriv_bmat3}. The same tensor ordering and Voigt convention are used for all generated material routines.
\begin{lstlisting}[
    style=fortranstyle,
    caption={Variable definition and calculation of deformation measures.},
    label={lst:pann_definition},
    firstnumber=107,
    breaklines=true,
]
  subroutine pann_eval_from_F(F, sig, K_out, G_out)
    implicit none
    real(dp), intent(in)  :: F(9)
    real(dp), intent(out) :: sig(6), K_out, G_out

    real(dp) :: identity(6)
    real(dp) :: b(6), b2(6)
    real(dp) :: I1, I2, I1b, I2b, sqrtI2b, J, Jm23, Jm43
    real(dp) :: dI1(6), dI1b(6), dI2b(6)
    real(dp) :: z1(N1), d1(N1)
    integer  :: m1, m2

    real(dp) :: g0(2)
    real(dp) :: g1(N1)

    call ensure_nn_params_loaded()
    K_out = K
    G_out = G

    ! b in Voigt: [b11,b22,b33,b12,b23,b13]
    b(1) = F(1)*F(1) + F(4)*F(4) + F(6)*F(6)
    b(2) = F(7)*F(7) + F(2)*F(2) + F(5)*F(5)
    b(3) = F(9)*F(9) + F(8)*F(8) + F(3)*F(3)
    b(4) = F(1)*F(7) + F(4)*F(2) + F(5)*F(6)
    b(5) = F(7)*F(9) + F(8)*F(2) + F(3)*F(5)
    b(6) = F(1)*F(9) + F(4)*F(8) + F(3)*F(6)

    ! b^2 in Voigt
    b2(1) = b(1)*b(1) + b(4)*b(4) + b(6)*b(6)
    b2(2) = b(4)*b(4) + b(2)*b(2) + b(5)*b(5)
    b2(3) = b(6)*b(6) + b(5)*b(5) + b(3)*b(3)
    b2(4) = b(1)*b(4) + b(2)*b(4) + b(5)*b(6)
    b2(5) = b(4)*b(6) + b(5)*b(2) + b(3)*b(5)
    b2(6) = b(1)*b(6) + b(4)*b(5) + b(3)*b(6)

    identity = [1.0_dp, 1.0_dp, 1.0_dp, 0.0_dp, 0.0_dp, 0.0_dp]

    I1  = b(1) + b(2) + b(3)
    dI1 = identity

    I2  = b(1)*b(2) + b(1)*b(3) + b(2)*b(3) - b(4)**2 - b(5)**2 - b(6)**2

    J = F(1)*F(2)*F(3) + F(4)*F(5)*F(9) + F(6)*F(7)*F(8) - F(9)*F(2)*F(6) - F(8)*F(5)*F(1) - F(3)*F(7)*F(4)
    Jm23 = J**(-2.0_dp/3.0_dp)
    Jm43 = Jm23*Jm23

    I1b  = Jm23 * I1
    dI1b = Jm23*(b-1.0_dp/3.0_dp*I1*identity)

    I2b  = Jm43 * I2
    sqrtI2b = sqrt(I2b)
    I2b = I2b * sqrtI2b
    dI2b = (1.5_dp * sqrtI2b) * Jm43*(I1*b - b2 - (2.0_dp/3.0_dp)*I2*identity)
\end{lstlisting}

Only quantities required for the final stress calculation are declared and evaluated. For example, the full derivative of $I_2$ is not explicitly formed. Instead, the modified second invariant $\overline{I}_2^\ast$ is computed as \texttt{I2b}, and the expression for its derivative already multiplied by $\mathbf{b}$ is evaluated as \texttt{dI2b}. The vector \texttt{identity} represents the second-order identity tensor in Voigt notation and is used wherever $\identMat$ appears in Eqs.~\eqref{eq:Deriv_bmat1}--\eqref{eq:Deriv_bmat3}.

The stress can now be evaluated according to Eq.~\eqref{eq:PK2_PANN}. The network derivative is obtained by a forward pass through the hidden layers followed by a backpropagation step to the invariant inputs. For the one-hidden-layer network considered here, only the derivative of the activation function is required in the final stress expression. The remaining material model reads:
\newpage
\begin{lstlisting}[
    style=fortranstyle,
    caption={Cauchy stress calculation for the PANN model.},
    label={lst:pann_stress},
    firstnumber=135,
    breaklines=true,
]
    ! --- NN forward pass ---
    ! Layer 1
    do m1=1,N1
      z1(m1) = b1_VEC(m1)
      z1(m1) = z1(m1) + W1_MAT(m1,1)*I1b + W1_MAT(m1,2)*I2b
      call softplusD(z1(m1), d1(m1))
    end do

    ! --- NN backward: compute dpsi_nn/dInp (size 2) ---
    ! Start with g1 = W2(1,:)
    g1 = W2_MAT(1, :)

    ! Backprop through layer 1 to input
    g0(1) = sum(g1*d1*W1_MAT(:,1))
    g0(2) = sum(g1*d1*W1_MAT(:,2))

    ! Svoigt from NN: 2 * [dpsi/dI1b, dpsi/dI2b] * [dI1b; dI2b]
    sig = 2.0_dp * ( g0(1)*dI1b + g0(2)*dI2b )

    ! Volume term for psi_vol = 0.5*K*(J-1)^2
    sig = sig/J + K*(J-1.0_dp)*dI1

  end subroutine pann_eval_from_F9

end module pann_model_001
\end{lstlisting}

With this procedure, the pretrained PANN is transferred into a standalone user material for \textsc{Radioss}. For comparison, a Carroll material model was implemented using the same tensor ordering and an analogous code structure. The full Fortran routine of the Carroll model is provided in Appendix~\ref{App:Subroutine_Carroll}.

\subsection{SoftPlus vs. SQuarePlus}

As discussed in Section~\ref{sec:PANNmod}, the SoftPlus function is commonly used in PANN architectures because it is smooth, convex, and monotonically increasing. However, its evaluation requires logarithmic and exponential functions, which can become expensive in a material routine that is called many times during an explicit simulation. Therefore, we additionally investigate the SQuarePlus function, introduced by Barron \cite{Barron.2021}, as an alternative activation function,
\revI{
\begin{align}
    \mathcal{SQP}(x)&=\cfrac{1}{2}\,\left(x+\sqrt{x^2+2}\right),\\
    \mathcal{SQP}'(x)&=\cfrac{1}{2}\,\left(1+\cfrac{x}{\sqrt{x^2+2}}\right).
\end{align}
}
The SQuarePlus function is also smooth, convex, and monotonically increasing, but replaces the logarithmic and exponential evaluations by an algebraic expression involving a square root. In the generated material routine, the network architecture and the constitutive formulation remain unchanged. The activation function used in the material routine is selected by an argument during the export of the Fortran user material and should correspond to the activation function used during training of the network. The routine \texttt{squareplus} evaluates both the SQuarePlus function and its derivative, whereas \texttt{squareplusD} evaluates only the derivative of the activation function. The corresponding implementations are shown in Listings~\ref{lst:pann_activation_SQP} and~\ref{lst:pann_activation_SQPD}.
\begin{lstlisting}[
    style=fortranstyle,
    caption={SQuarePlu activation function and derivative of the generated Fortran material module.},
    label={lst:pann_activation_SQP},
    firstnumber=42,
]
  pure elemental subroutine squareplus(x, y, dy)
    real(dp), intent(in)  :: x
    real(dp), intent(out) :: y, dy
    real(dp), parameter   :: b = 2.0_dp
    real(dp) :: s
      s = sqrt(x*x + b)
      y  = 0.5_dp * (x + s)
      dy = 0.5_dp * (1.0_dp + x / s)
  end subroutine squareplus
\end{lstlisting}
\begin{lstlisting}[
    style=fortranstyle,
    caption={Derivative-only SQuarePlu activation routine of the generated Fortran material module.},
    label={lst:pann_activation_SQPD},
    firstnumber=52,
]
  pure elemental subroutine squareplusD(x, dy)
    real(dp), intent(in)  :: x
    real(dp), intent(out) :: dy
    real(dp), parameter   :: b = 2.0_dp
    real(dp) :: s
      s = sqrt(x*x + b)
      dy = 0.5_dp * (1.0_dp + x / s)
  end subroutine squareplusD
\end{lstlisting}

To quantify the isolated computational difference between the SoftPlus and SQuarePlus activation functions, a small benchmark test was performed. Both activation functions and their derivatives were evaluated repeatedly inside a loop, and the resulting values were accumulated to prevent compiler optimization from removing the function calls. In addition, a baseline measurement without any activation-function evaluation was carried out to account for the loop overhead. The full Fortran benchmark code is provided in Appendix~\ref{App:Benchmark_SoftSquare}. On our system, the net evaluation time of the SQuarePlus routine was reduced by approximately $87$--$88\,\%$ compared with the SoftPlus routine after subtracting the baseline runtime. The influence of this replacement on the complete finite element simulation is discussed in the following section.

For the runtime comparison presented in Section~5.2, the SQuarePlus-based PANNs were trained separately using the same training procedure and network architectures as the corresponding SoftPlus-based models. Thus, the reported SQuarePlus results do not rely on a post-processing replacement of the activation function after training, but on networks optimized directly with the SQuarePlus activation. The resulting fitted material responses were visually indistinguishable from the SoftPlus-based fits for the considered deformation modes. For completeness, we also tested replacing the activation function only after training. For the shallow network with one hidden layer, applying the SQuarePlus function to a network originally trained with SoftPlus led to very similar material responses. However, for deeper architectures the deviations became more pronounced, since small differences introduced by the changed activation function are propagated through multiple nonlinear layers. Therefore, all SQuarePlus models used in the following runtime evaluation were trained directly with the SQuarePlus activation.
\section{Computational example and runtime evaluation}
\label{sec:Results}
Based on the model definition and computational setup described in Section~\ref{sec:Setup_definition}, the impact of the rigid projectile on the dragon model is simulated using the different material formulations. For each of the three mesh configurations, all material formulations are simulated with identical boundary conditions, mesh configurations, hardware settings, and solver settings. The runtime values reported below are averaged over four identical runs in order to reduce the influence of random timing variations. In the following, the simulation results are first compared qualitatively before the computational efficiency of the material formulations is evaluated. Snapshots of the simulation sequence are provided in Appendix~\ref{App:Sequence}.

The different PANN architectures are denoted by expressions such as $\mathrm{PANN}_{\mathcal{SP}}\:(5)$ and $\mathrm{PANN}_{\mathcal{SQP}}\:(5-5)$. Here, the subscript identifies the activation function, with $\mathcal{SP}$ and $\mathcal{SQP}$ referring to SoftPlus and SQuarePlus, respectively, while the numbers in parentheses specify the hidden-layer architecture.

\subsection{Qualitative comparison of the stress response}
\label{sec:Comparison}
The PANN models and the Carroll model do not yield identical fitted material responses and also differ slightly in their initial shear moduli. Therefore, exact agreement of the stress and deformation fields is not expected, particularly in local point-wise quantities. Instead, the simulations are used to assess whether the generated PANN material routines achieve numerically stable simulations and produce mechanically consistent stress and displacement fields that are comparable to the classical reference model.

The PANN models and the Carroll model do not yield identical fitted material responses and also differ slightly in their initial shear moduli. Therefore, exact agreement of the stress and deformation fields is not expected, particularly in local point-wise quantities. Instead, the simulations are used to assess whether the generated PANN material routines lead to numerically stable simulations and produce mechanically consistent stress fields and deformation patterns that are comparable to those obtained with the classical reference model.

For this purpose, the von Mises stress distribution is evaluated at $t_1=2\,\mathrm{ms}$, after the projectile has rebounded and a pronounced stress state has developed in the dragon head. For better visualization, the color scale is limited to $0.5\,\mathrm{MPa}$, although the maximum von Mises stress $\sigma_{\mathrm{vM},\max}$ is higher in all simulations. The results are shown for the Carroll model, a PANN with one hidden layer and five neurons ($\mathrm{PANN}_{\mathcal{SP}}\:(5)$), and a PANN with two hidden layers and five neurons each ($\mathrm{PANN}_{\mathcal{SP}}\:(5-5)$). The corresponding results for the three mesh resolutions are shown in Figs.~\ref{fig:mesh1_von_mises}--\ref{fig:mesh3_von_mises}. The contour plots show only the results of the networks trained with SoftPlus activation function, since no visible differences were observed in the corresponding networks trained with SQuarePlus.
\begin{figure}[htb]
    \centering
    \begin{minipage}[t]{0.27\textwidth}
        \centering
        {\small Carroll}
        \vspace{1mm}
        \includegraphics[
            width=\linewidth
        ]{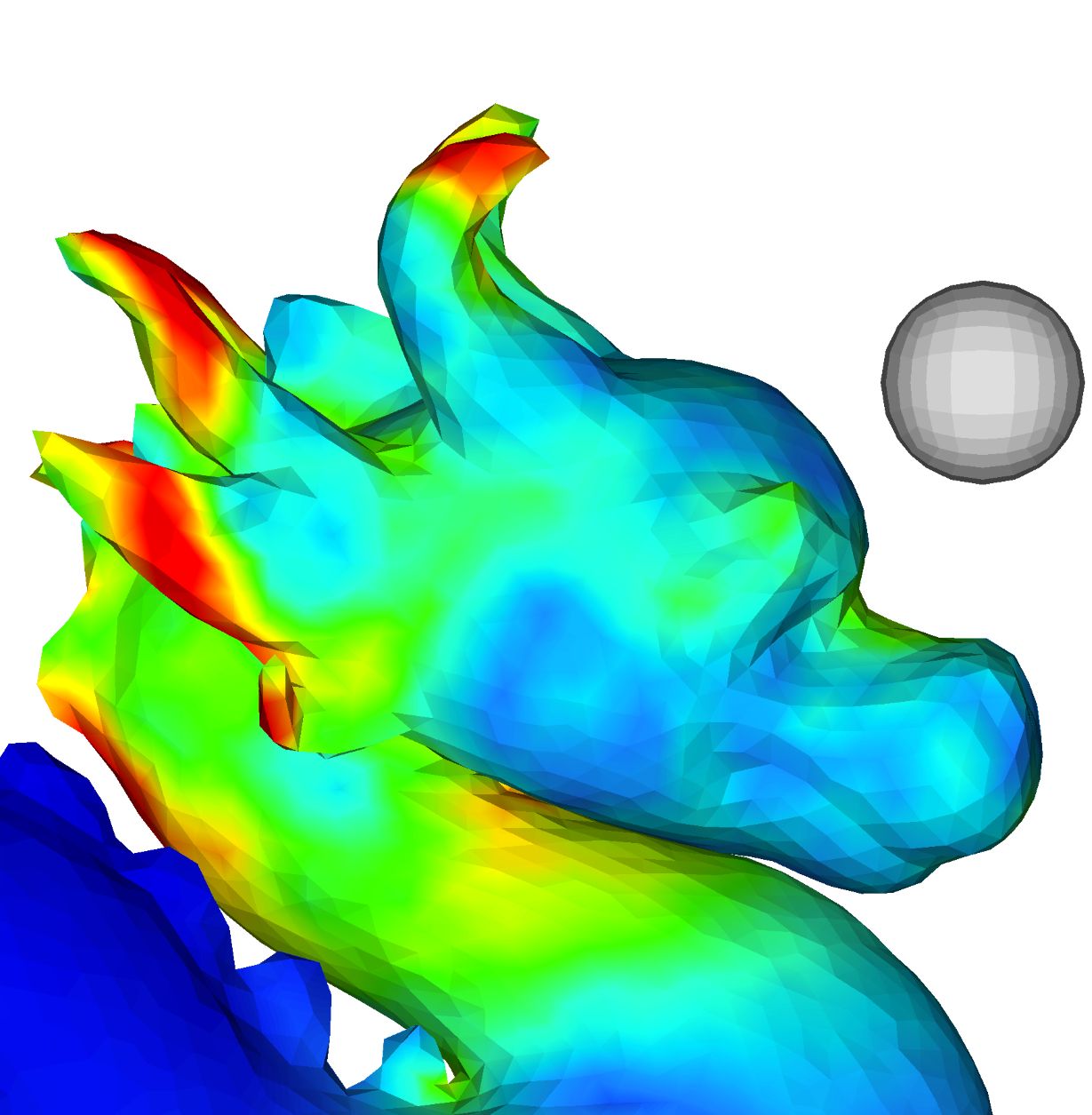}
        \vspace{1mm}
        $\sigma_{\mathrm{vM},\max}=1.706 \,\mathrm{MPa}$
    \end{minipage}
    \hfill
    \begin{minipage}[t]{0.27\textwidth}
        \centering
        {\small $\mathrm{PANN}_{\mathcal{SP}}\:(5)$}
        \vspace{1mm}
        \includegraphics[
            width=\linewidth
        ]{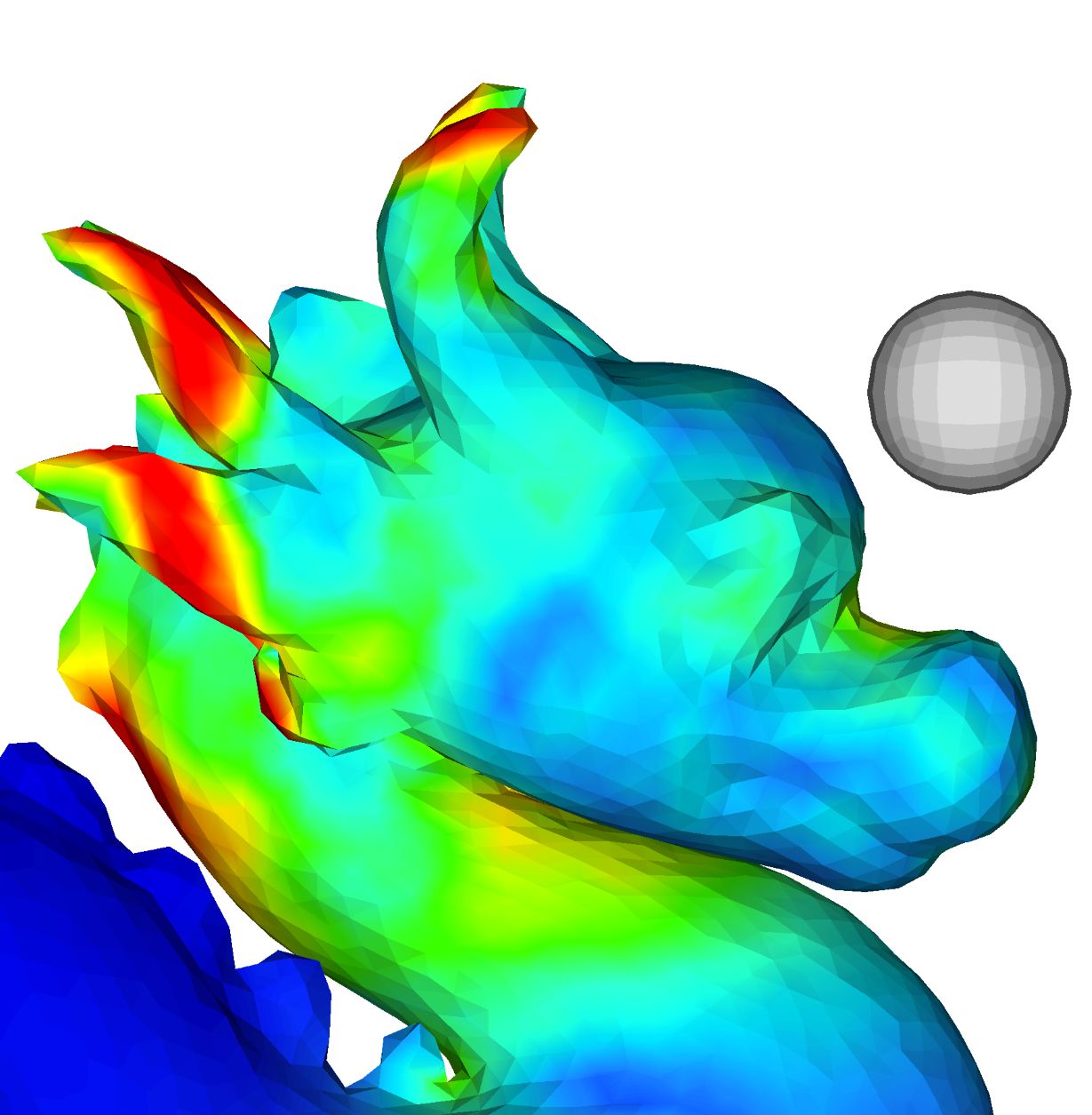}
        \vspace{1mm}
        $\sigma_{\mathrm{vM},\max}=1.718 \,\mathrm{MPa}$
    \end{minipage}
    \hfill
    \begin{minipage}[t]{0.27\textwidth}
        \centering
        {\small $\mathrm{PANN}_{\mathcal{SP}}\:(5-5)$}
        \vspace{1mm}
        \includegraphics[
            width=\linewidth
        ]{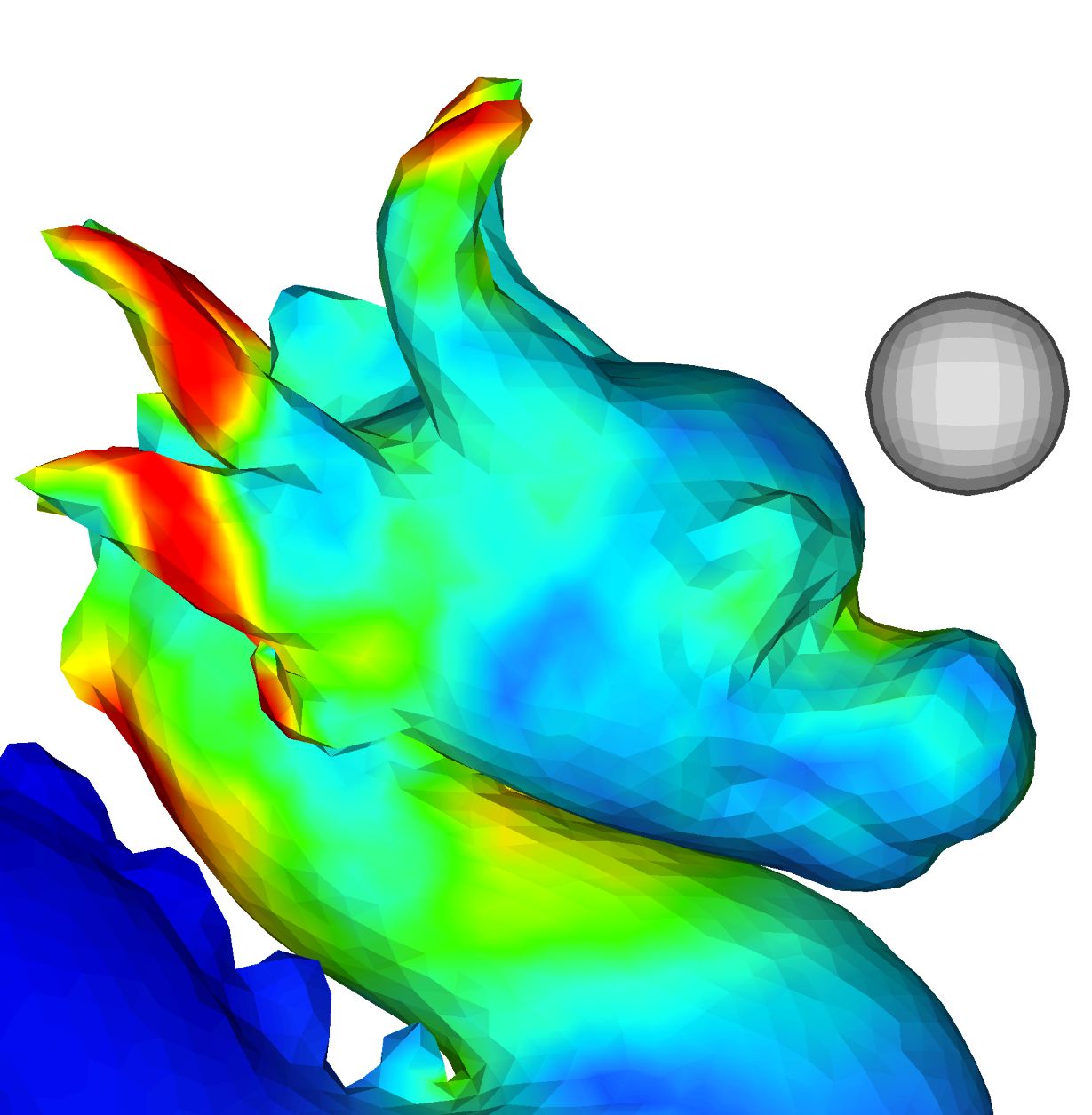}
        \vspace{1mm}
        $\sigma_{\mathrm{vM},\max}=1.741 \,\mathrm{MPa}$
    \end{minipage}
    \hfill
    \begin{minipage}[t]{0.13\textwidth}
        \centering
        \vspace{-10mm} % anpassen
        \setlength{\tkwi}{0.9\textwidth}
        \setlength{\tkhe}{2.5\tkwi}
        \definecolor{c1}{rgb}{1.00000,0.00000,0.00000}%
\definecolor{c2}{rgb}{1.00000,0.84300,0.00000}%
\definecolor{c3}{rgb}{0.44700,1.00000,0.00000}%
\definecolor{c4}{rgb}{0.19600,1.00000,0.41200}%
\definecolor{c5}{rgb}{0.01600,0.74100,0.91400}%
\definecolor{c6}{rgb}{0.00400,0.02700,0.84300}%

\begin{tikzpicture}

% Zeichenfläche für Skalierung über \tkwi und \tkhe
\path[use as bounding box]
    (0\tkwi,0\tkhe) rectangle (1\tkwi,1\tkhe);

% Titel
%\node[
%    anchor=north west,
%    align=left
%] at (0.00\tkwi,1.00\tkhe)
%{Von Mises Stress\\Simple Average};

% Segment 1
\shade[
    top color=c1,
    bottom color=c2,
    draw=none
] (0.00\tkwi,0.78\tkhe) rectangle (0.18\tkwi,0.624\tkhe);

% Segment 2
\shade[
    top color=c2,
    bottom color=c3,
    draw=none
] (0.00\tkwi,0.624\tkhe) rectangle (0.18\tkwi,0.469\tkhe);

% Segment 3
\shade[
    top color=c3,
    bottom color=c4,
    draw=none
] (0.00\tkwi,0.469\tkhe) rectangle (0.18\tkwi,0.391\tkhe);

% Segment 4
\shade[
    top color=c4,
    bottom color=c5,
    draw=none
] (0.00\tkwi,0.391\tkhe) rectangle (0.18\tkwi,0.2356\tkhe);

% Segment 5
\shade[
    top color=c5,
    bottom color=c6,
    draw=none
] (0.00\tkwi,0.2356\tkhe) rectangle (0.18\tkwi,0.08\tkhe);

% Maximalwert oben
\node[
    anchor=west
] at (0.26\tkwi,0.78\tkhe)
{$0.5\,\mathrm{MPa}$};

% Minimalwert unten
\node[
    anchor=west
] at (0.26\tkwi,0.08\tkhe)
{$0.0\,\mathrm{MPa}$};

\end{tikzpicture}
    \end{minipage}
    \caption{Von Mises stress distribution for Mesh 1 with a characteristic edge length of $1.6\,\mathrm{mm}$ using the Carroll model as well as the $\mathrm{PANN}_{\mathcal{SP}}\:(5)$ and $\mathrm{PANN}_{\mathcal{SP}}\:(5-5)$ models.}
    \label{fig:mesh1_von_mises}
\end{figure}
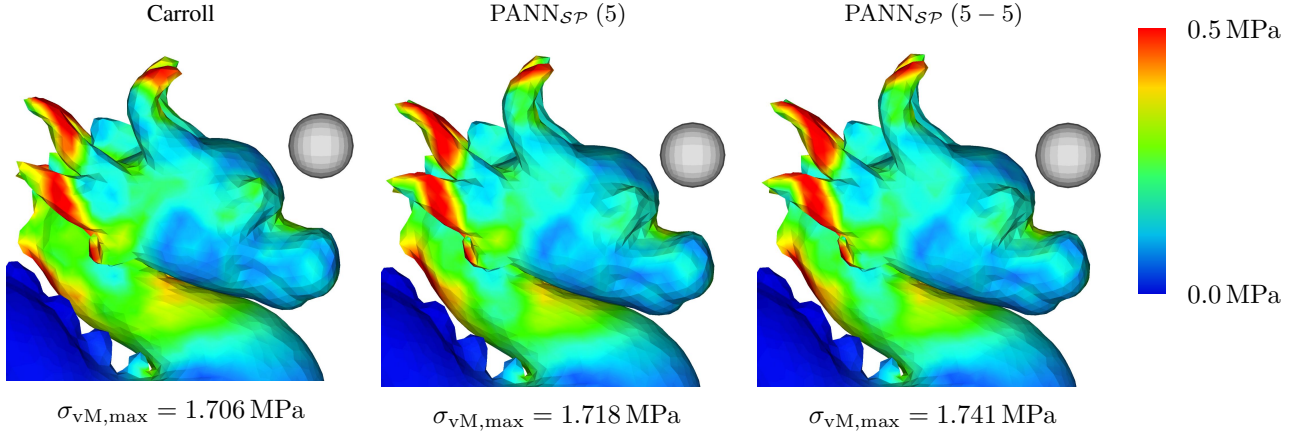
\begin{figure}[htb]
    \centering
    \begin{minipage}[t]{0.27\textwidth}
        \centering
        {\small Carroll}
        \vspace{1mm}
        \includegraphics[
            width=\linewidth
        ]{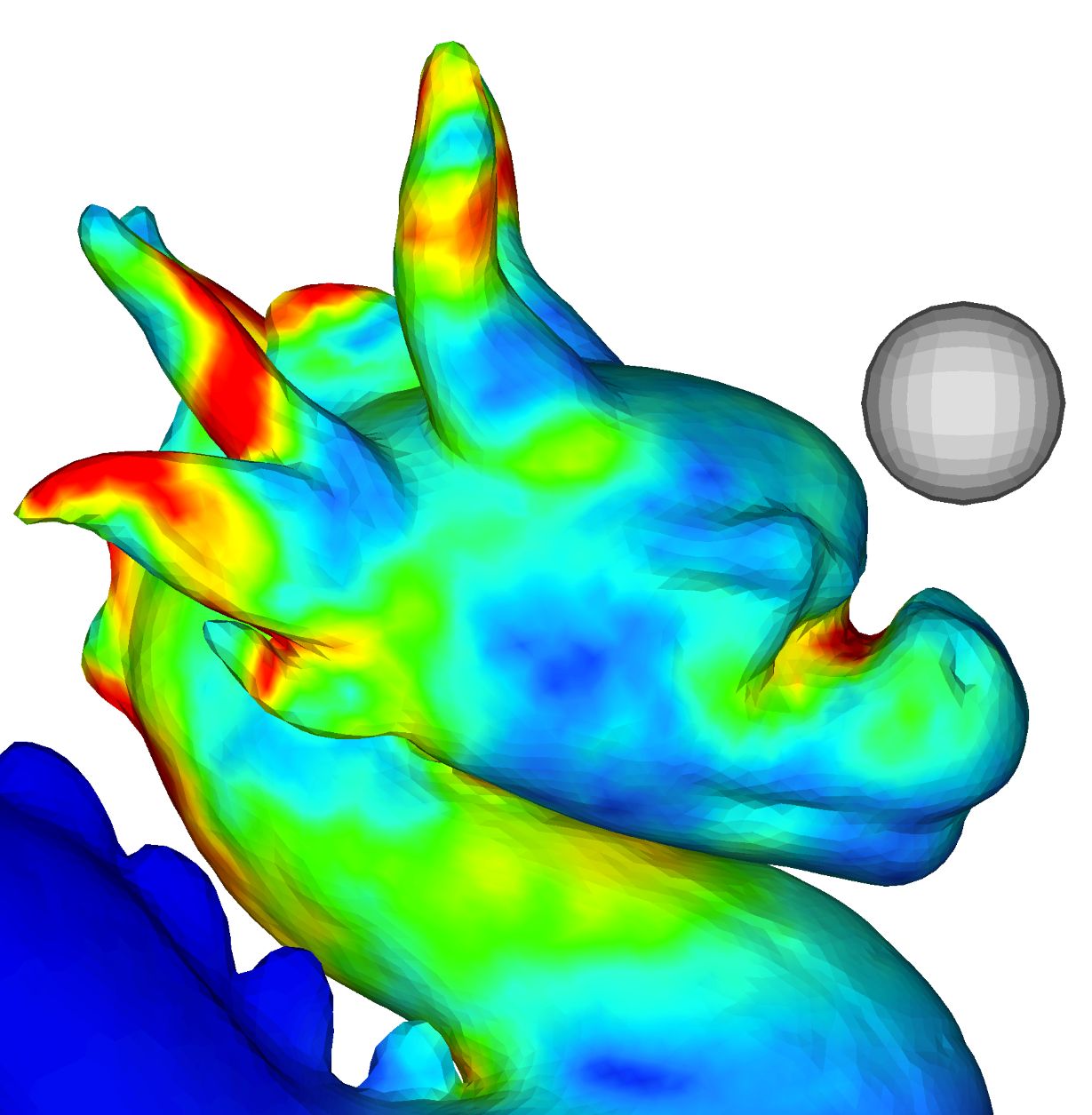}
        \vspace{1mm}
        $\sigma_{\mathrm{vM},\max}=3.383 \,\mathrm{MPa}$
    \end{minipage}
    \hfill
    \begin{minipage}[t]{0.27\textwidth}
        \centering
        {\small $\mathrm{PANN}_{\mathcal{SP}}\:(5)$}
        \vspace{1mm}
        \includegraphics[
            width=\linewidth
        ]{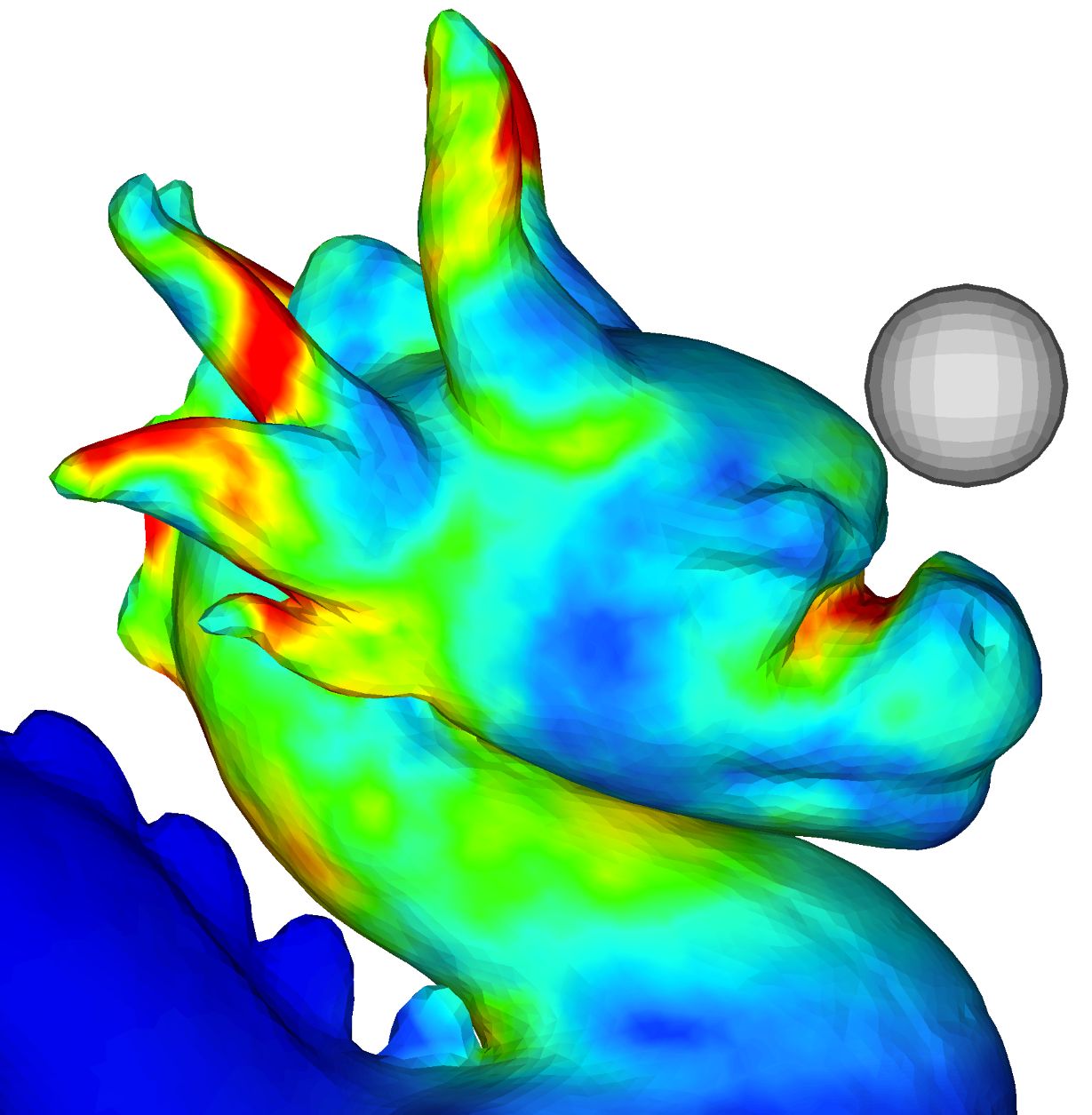}
        \vspace{1mm}
        $\sigma_{\mathrm{vM},\max}=3.745 \,\mathrm{MPa}$
    \end{minipage}
    \hfill
    \begin{minipage}[t]{0.27\textwidth}
        \centering
        {\small $\mathrm{PANN}_{\mathcal{SP}}\:(5-5)$}
        \vspace{1mm}
        \includegraphics[
            width=\linewidth
        ]{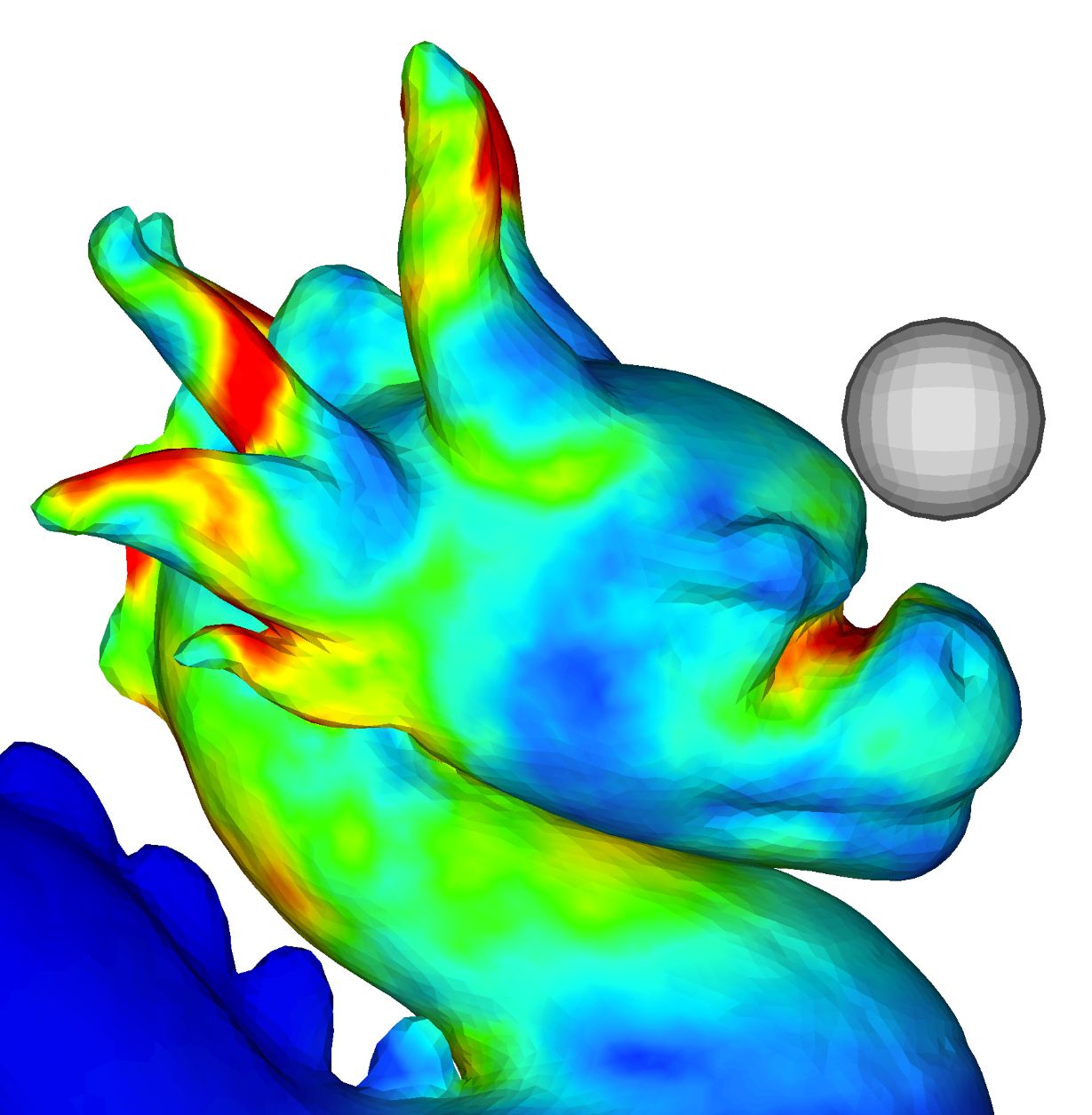}
        \vspace{1mm}
        $\sigma_{\mathrm{vM},\max}=3.769 \,\mathrm{MPa}$
    \end{minipage}
    \hfill
    \begin{minipage}[t]{0.13\textwidth}
        \centering
        \vspace{-10mm} % anpassen
        \setlength{\tkwi}{0.9\textwidth}
        \setlength{\tkhe}{2.5\tkwi}
        \definecolor{c1}{rgb}{1.00000,0.00000,0.00000}%
\definecolor{c2}{rgb}{1.00000,0.84300,0.00000}%
\definecolor{c3}{rgb}{0.44700,1.00000,0.00000}%
\definecolor{c4}{rgb}{0.19600,1.00000,0.41200}%
\definecolor{c5}{rgb}{0.01600,0.74100,0.91400}%
\definecolor{c6}{rgb}{0.00400,0.02700,0.84300}%

\begin{tikzpicture}

% Zeichenfläche für Skalierung über \tkwi und \tkhe
\path[use as bounding box]
    (0\tkwi,0\tkhe) rectangle (1\tkwi,1\tkhe);

% Titel
%\node[
%    anchor=north west,
%    align=left
%] at (0.00\tkwi,1.00\tkhe)
%{Von Mises Stress\\Simple Average};

% Segment 1
\shade[
    top color=c1,
    bottom color=c2,
    draw=none
] (0.00\tkwi,0.78\tkhe) rectangle (0.18\tkwi,0.624\tkhe);

% Segment 2
\shade[
    top color=c2,
    bottom color=c3,
    draw=none
] (0.00\tkwi,0.624\tkhe) rectangle (0.18\tkwi,0.469\tkhe);

% Segment 3
\shade[
    top color=c3,
    bottom color=c4,
    draw=none
] (0.00\tkwi,0.469\tkhe) rectangle (0.18\tkwi,0.391\tkhe);

% Segment 4
\shade[
    top color=c4,
    bottom color=c5,
    draw=none
] (0.00\tkwi,0.391\tkhe) rectangle (0.18\tkwi,0.2356\tkhe);

% Segment 5
\shade[
    top color=c5,
    bottom color=c6,
    draw=none
] (0.00\tkwi,0.2356\tkhe) rectangle (0.18\tkwi,0.08\tkhe);

% Maximalwert oben
\node[
    anchor=west
] at (0.26\tkwi,0.78\tkhe)
{$0.5\,\mathrm{MPa}$};

% Minimalwert unten
\node[
    anchor=west
] at (0.26\tkwi,0.08\tkhe)
{$0.0\,\mathrm{MPa}$};

\end{tikzpicture}
    \end{minipage}
    \caption{Von Mises stress distribution for Mesh 2 with a characteristic edge length of $0.8\,\mathrm{mm}$ using the Carroll model as well as the $\mathrm{PANN}_{\mathcal{SP}}\:(5)$ and $\mathrm{PANN}_{\mathcal{SP}}\:(5-5)$ models.}
    \label{fig:mesh2_von_mises}
\end{figure}
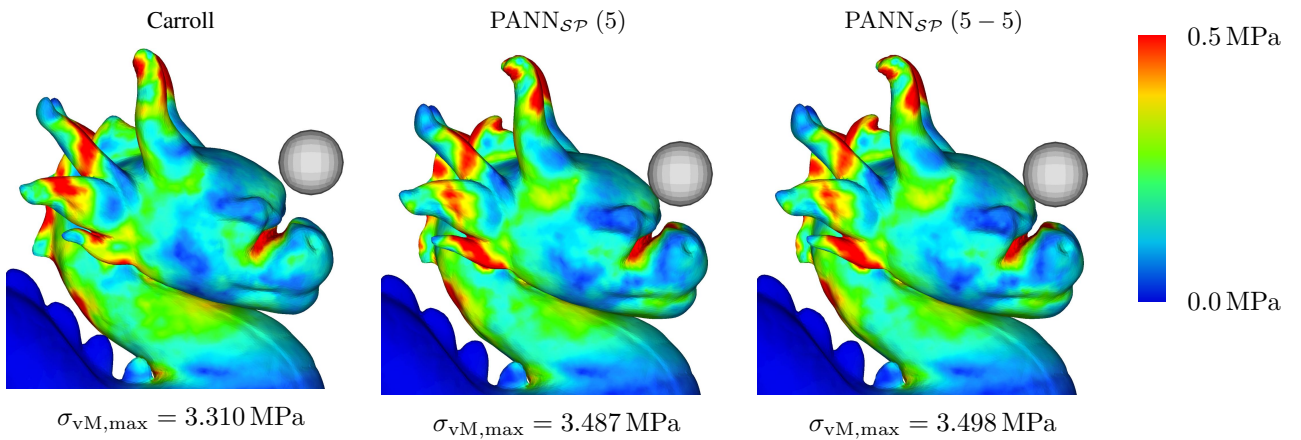
\begin{figure}[H]
    \centering
    \begin{minipage}[t]{0.27\textwidth}
        \centering
        {\small Carroll}
        \vspace{1mm}
        \includegraphics[
            width=\linewidth
        ]{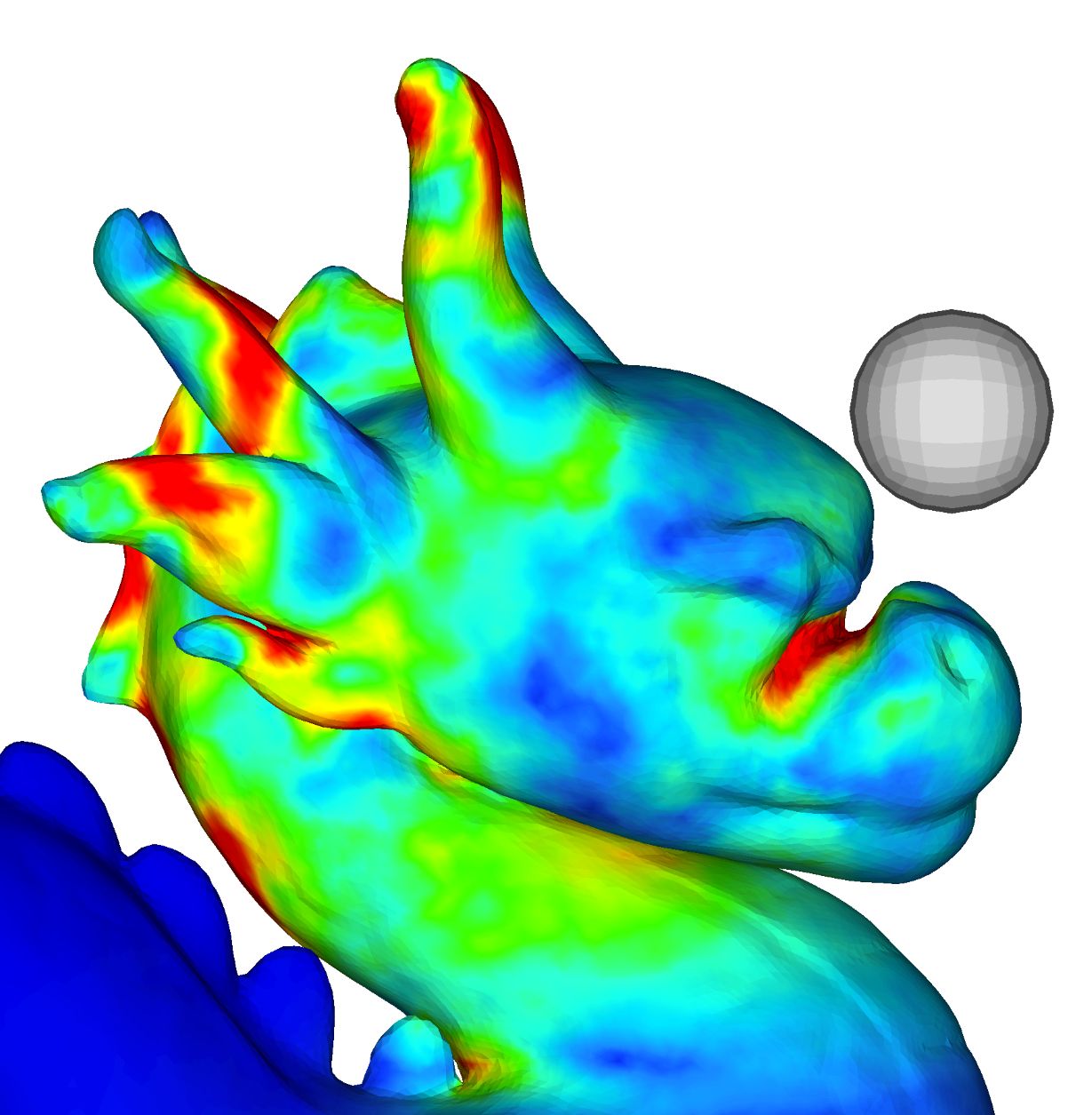}
        \vspace{1mm}
        $\sigma_{\mathrm{vM},\max}=3.310 \,\mathrm{MPa}$
    \end{minipage}
    \hfill
    \begin{minipage}[t]{0.27\textwidth}
        \centering
        {\small $\mathrm{PANN}_{\mathcal{SP}}\:(5)$}
        \vspace{1mm}
        \includegraphics[
            width=\linewidth
        ]{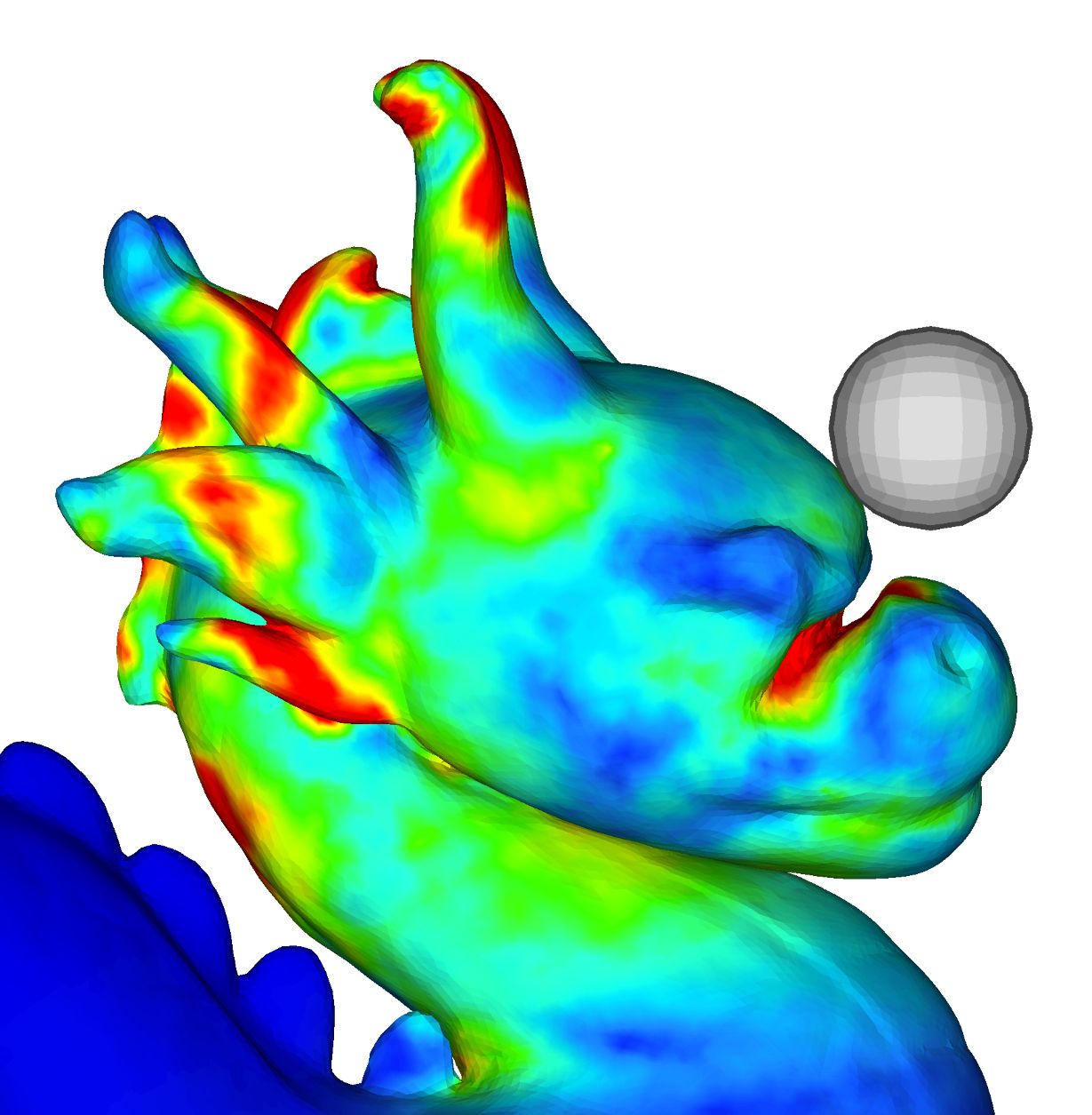}
        \vspace{1mm}
        $\sigma_{\mathrm{vM},\max}=3.487 \,\mathrm{MPa}$
    \end{minipage}
    \hfill
    \begin{minipage}[t]{0.27\textwidth}
        \centering
        {\small $\mathrm{PANN}_{\mathcal{SP}}\:(5-5)$}
        \vspace{1mm}
        \includegraphics[
            width=\linewidth
        ]{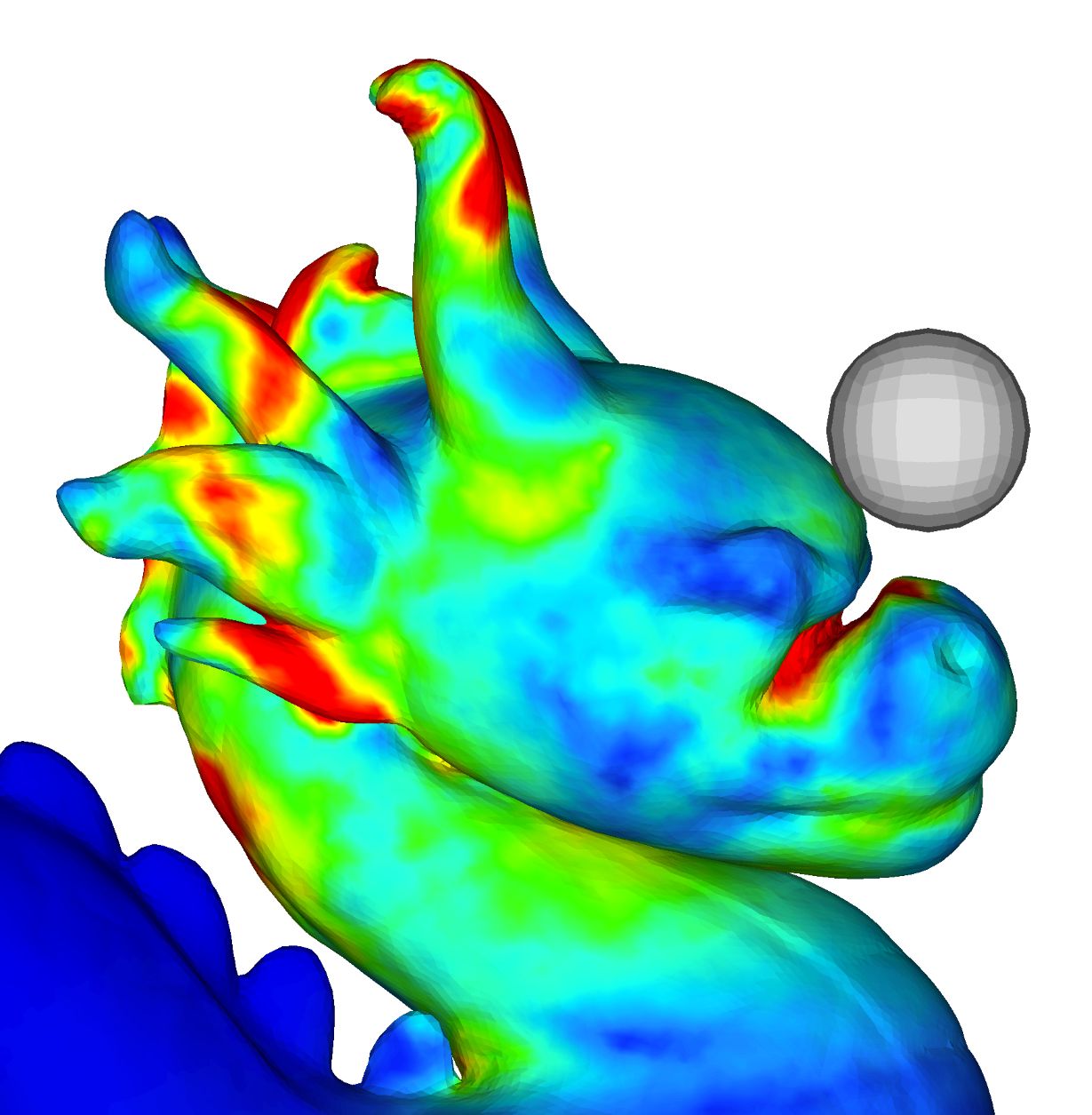}
        \vspace{1mm}
        $\sigma_{\mathrm{vM},\max}=3.498 \,\mathrm{MPa}$
    \end{minipage}
    \hfill
    \begin{minipage}[t]{0.13\textwidth}
        \centering
        \vspace{-10mm} % anpassen
        \setlength{\tkwi}{0.9\textwidth}
        \setlength{\tkhe}{2.5\tkwi}
        \definecolor{c1}{rgb}{1.00000,0.00000,0.00000}%
\definecolor{c2}{rgb}{1.00000,0.84300,0.00000}%
\definecolor{c3}{rgb}{0.44700,1.00000,0.00000}%
\definecolor{c4}{rgb}{0.19600,1.00000,0.41200}%
\definecolor{c5}{rgb}{0.01600,0.74100,0.91400}%
\definecolor{c6}{rgb}{0.00400,0.02700,0.84300}%

\begin{tikzpicture}

% Zeichenfläche für Skalierung über \tkwi und \tkhe
\path[use as bounding box]
    (0\tkwi,0\tkhe) rectangle (1\tkwi,1\tkhe);

% Titel
%\node[
%    anchor=north west,
%    align=left
%] at (0.00\tkwi,1.00\tkhe)
%{Von Mises Stress\\Simple Average};

% Segment 1
\shade[
    top color=c1,
    bottom color=c2,
    draw=none
] (0.00\tkwi,0.78\tkhe) rectangle (0.18\tkwi,0.624\tkhe);

% Segment 2
\shade[
    top color=c2,
    bottom color=c3,
    draw=none
] (0.00\tkwi,0.624\tkhe) rectangle (0.18\tkwi,0.469\tkhe);

% Segment 3
\shade[
    top color=c3,
    bottom color=c4,
    draw=none
] (0.00\tkwi,0.469\tkhe) rectangle (0.18\tkwi,0.391\tkhe);

% Segment 4
\shade[
    top color=c4,
    bottom color=c5,
    draw=none
] (0.00\tkwi,0.391\tkhe) rectangle (0.18\tkwi,0.2356\tkhe);

% Segment 5
\shade[
    top color=c5,
    bottom color=c6,
    draw=none
] (0.00\tkwi,0.2356\tkhe) rectangle (0.18\tkwi,0.08\tkhe);

% Maximalwert oben
\node[
    anchor=west
] at (0.26\tkwi,0.78\tkhe)
{$0.5\,\mathrm{MPa}$};

% Minimalwert unten
\node[
    anchor=west
] at (0.26\tkwi,0.08\tkhe)
{$0.0\,\mathrm{MPa}$};

\end{tikzpicture}
    \end{minipage}
    \caption{Von Mises stress distribution for Mesh 3 with a characteristic edge length of $0.5\,\mathrm{mm}$ using the Carroll model as well as the $\mathrm{PANN}_{\mathcal{SP}}\:(5)$ and $\mathrm{PANN}_{\mathcal{SP}}\:(5-5)$ models.}
    \label{fig:mesh3_von_mises}
\end{figure}
Overall, the stress distributions obtained with the PANN-based material models are qualitatively close to those obtained with the Carroll model. The coarsest mesh shows a visibly different deformation pattern and lower maximum stresses than the two finer discretizations. The two finer meshes still exhibit local deviations, but the deformation shape and stress distribution are similar. A more detailed mesh convergence study would be useful, but is beyond the scope of the present work. For the purpose of this study, the results show that the generated PANN user materials can be applied to the explicit impact simulation without numerical instabilities and yield physically plausible stress fields.

\begin{figure}[htb]
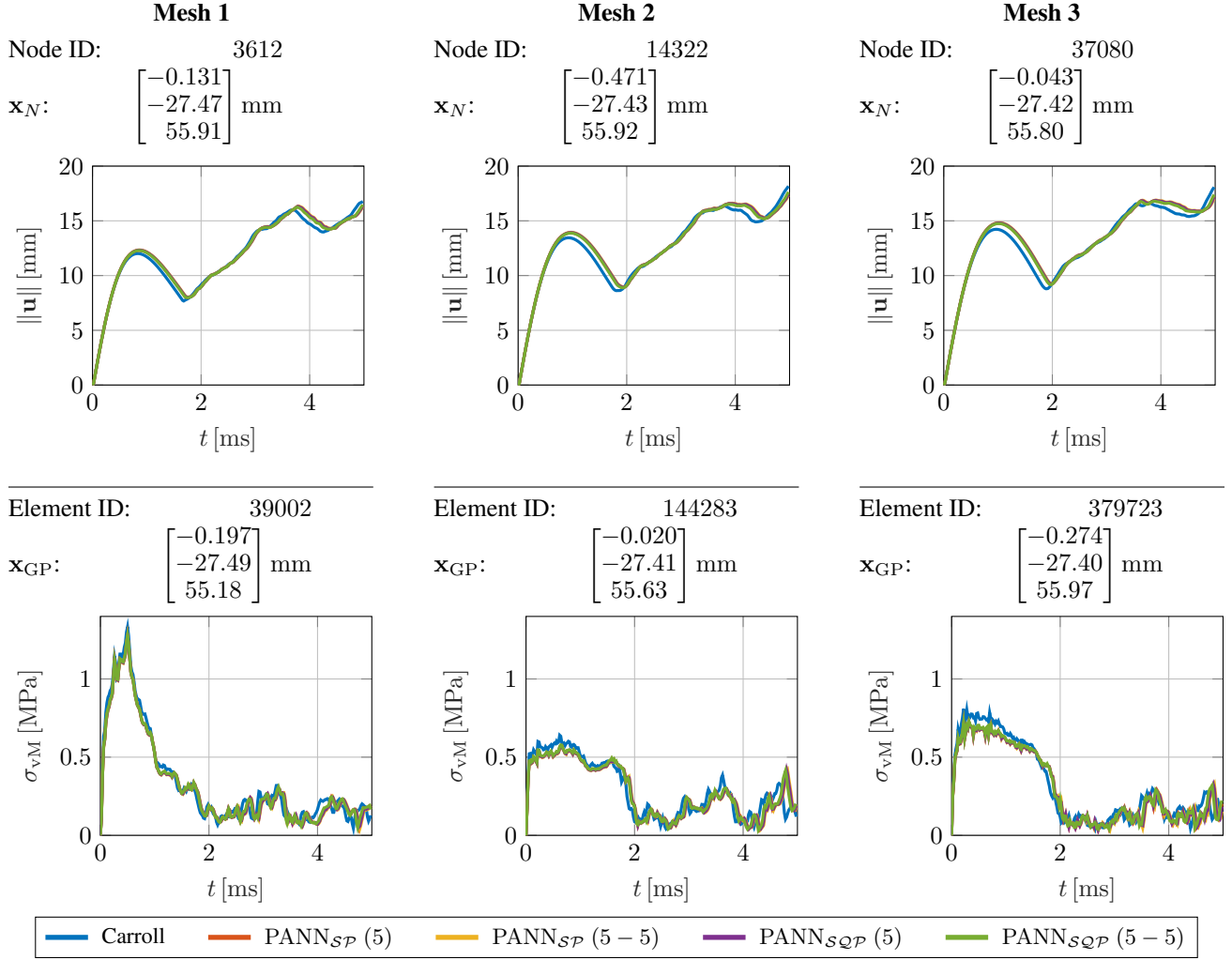

    \centering
    \begin{minipage}[t]{0.3\textwidth}
        \centering
        \textbf{Mesh 1}\\
        \vspace{2mm}
        {\raggedright
        \begin{tabular}{@{}l r@{}}
        Node ID: & $3612$\\
        $\mathbf{x}_N$: & $\begin{bmatrix}
        -0.131\\-27.47 \\ \phantom{-}55.91
        \end{bmatrix}\,\mathrm{mm}$
        \end{tabular}
        \par}
        \setlength\tkwi{0.75\textwidth}
        \setlength\tkhe{0.8\tkwi}
        \input{tikz/Mesh08_Displ_N3612.tikz}
        \vspace{1mm}
        \noindent\rule{\textwidth}{0.4pt}
        \vspace{1mm}
        {\raggedright
        \begin{tabular}{@{}l r@{}}
        Element ID: & $39002$\\
        $\mathbf{x}_{\mathrm{GP}}$: & $\begin{bmatrix}
        -0.197\\-27.49 \\  55.18
        \end{bmatrix}\,\mathrm{mm}$
        \end{tabular}
        \par}
        \input{tikz/Mesh08_Stress_E39002.tikz}
    \end{minipage}
    \hfill
    \begin{minipage}[t]{0.3\textwidth}
        \centering
        \textbf{Mesh 2}\\
        \vspace{2mm}
        {\raggedright
        \begin{tabular}{@{}l r@{}}
        Node ID: & $14322$\\
        $\mathbf{x}_N$: & $\begin{bmatrix}
        -0.471\\-27.43 \\ 55.92
        \end{bmatrix}\,\mathrm{mm}$
        \end{tabular}
        \par}
        \setlength\tkwi{0.75\textwidth}
        \setlength\tkhe{0.8\tkwi}
        \input{tikz/Mesh04_Displ_N14322.tikz}
        \vspace{1mm}
        \noindent\rule{\textwidth}{0.4pt}
        \vspace{1mm}
        {\raggedright
        \begin{tabular}{@{}l r@{}}
        Element ID: & $144283$\\
        $\mathbf{x}_{\mathrm{GP}}$: & $\begin{bmatrix}
        -0.020\\-27.41 \\  55.63
        \end{bmatrix}\,\mathrm{mm}$
        \end{tabular}
        \par}
        \input{tikz/Mesh04_Stress_E144283.tikz}
    \end{minipage}
    \hfill
    \begin{minipage}[t]{0.3\textwidth}
        \centering
        \textbf{Mesh 3}\\
        \vspace{2mm}
        {\raggedright
        \begin{tabular}{@{}l r@{}}
        Node ID: & $37080$\\
        $\mathbf{x}_N$: & $\begin{bmatrix}
        -0.043\\-27.42 \\  55.80
        \end{bmatrix}\,\mathrm{mm}$
        \end{tabular}
        \par}
        \setlength\tkwi{0.75\textwidth}
        \setlength\tkhe{0.8\tkwi}
        \input{tikz/Mesh025_Displ_N37080.tikz}
        \vspace{1mm}
        \noindent\rule{\textwidth}{0.4pt}
        \vspace{1mm}
        {\raggedright
        \begin{tabular}{@{}l r@{}}
        Element ID: & $379723$\\
        $\mathbf{x}_{\mathrm{GP}}$: & $\begin{bmatrix}
        -0.274\\-27.40 \\  55.97
        \end{bmatrix}\,\mathrm{mm}$
        \end{tabular}
        \par}
        \input{tikz/Mesh025_Stress_E379723.tikz}
    \end{minipage}
    \\[0.05cm]
    \centering
\definecolor{mycolor1}{rgb}{0.00000,0.44700,0.74100}%
\definecolor{mycolor2}{rgb}{0.85000,0.32500,0.09800}%
\definecolor{mycolor3}{rgb}{0.92900,0.69400,0.12500}%
\definecolor{mycolor4}{rgb}{0.49400,0.18400,0.55600}%
\definecolor{mycolor5}{rgb}{0.46600,0.67400,0.18800}%
\begin{tikzpicture}

\begin{axis}[
    hide axis,
    xmin=0, xmax=1,
    ymin=0, ymax=1,
    legend columns=5,               % <-- Zwei Spalten
    legend style={
        column sep=2.5pt,
        /tikz/every even column/.append style={column sep=0.5cm},
        font=\footnotesize,
        cells={anchor=west}
    }
]
\addlegendimage{color=mycolor1, line width=2.0pt}
\addlegendentry{Carroll}

\addlegendimage{color=mycolor2, line width=2.0pt}
\addlegendentry{$\mathrm{PANN}_{\mathcal{SP}}\:(5)$}

\addlegendimage{color=mycolor3, line width=2.0pt}
\addlegendentry{$\mathrm{PANN}_{\mathcal{SP}}\:(5-5)$}

\addlegendimage{color=mycolor4, line width=2.0pt}
\addlegendentry{$\mathrm{PANN}_{\mathcal{SQP}}\:(5)$}

\addlegendimage{color=mycolor5, line width=2.0pt}
\addlegendentry{$\mathrm{PANN}_{\mathcal{SQP}}\:(5-5)$}

\end{axis}
\end{tikzpicture}%
    \caption{Comparison of the five different material formulations for all three mesh discretizations. For each mesh, the displacement magnitude $\|\mathbf{u}\|$ at a selected node and the von Mises stress $\sigma_\mathrm{vM}$ at a selected element are shown over time. Differences between the mesh levels are attributed to the slightly different spatial locations of the selected nodes and elements in the respective discretizations.}
    \label{fig:Comparison_PANNs}
\end{figure}

The transient response is additionally compared in Fig.~\ref{fig:Comparison_PANNs} using the displacement magnitude at a selected node and the von Mises stress at a selected element for each mesh. The selected node and element correspond to the region that first comes into contact with the projectile. Since the discretizations differ between the meshes, the node and element were chosen at comparable spatial locations as closely as possible, although an exact correspondence is not possible. This comparison is therefore not intended as a strict point-wise evaluation, but as an additional assessment of the temporal response. The four PANN variants show very good agreement, while the Carroll response exhibits only minor deviations, which are expected due to the different fit to the experimental data discussed above.

\subsection{Runtime evaluation}
\label{sec:Runtime}
A central objective of the proposed implementation is the computational efficiency of the generated material routines. Data-driven material models can only be used in routine engineering simulations if their additional cost remains comparable to established constitutive models.

\textsc{Radioss} reports the CPU time separated into several contributions, including contact sorting, contact forces, element forces, kinematic conditions, integration, assembling, and other operations such as input/output (I/O). In addition, the number of explicit cycles $n_\mathrm{cy}$ is reported. For the runtime comparison, the finest mesh is considered, since it leads to the largest overall computational effort and is therefore less affected by small timing artifacts.

The most relevant quantity for the present comparison is the CPU time spent in the \emph{element forces} routine, denoted by $T_\mathrm{EF}$. This part contains the evaluation of the constitutive model and therefore reflects the additional cost of the PANN material formulation most directly. Since the number of cycles differs slightly between the simulations, the element-force time is normalized by the number of cycles according to
\begin{equation}
    t_\mathrm{EF}
    =
    \cfrac{T_\mathrm{EF}}{n_\mathrm{cy}} .
\end{equation}
The resulting normalized element-force times are summarized in Table~\ref{tab:runtime_element_forces}. The values represent the mean of four identical runs. The deviation between the individual runs was below $0.5\,\%$ in all cases.
\begin{table}[htbp]
    \centering
    \caption{Runtime comparison based on the CPU time spent in the element-force calculation for the finest mesh. The normalized time $t_\mathrm{EF}$ is obtained by dividing the element-force time by the number of explicit cycles.}
    \label{tab:runtime_element_forces}
    \begin{tabular}{lcccccc}
        \hline
        Model & Activation & Architecture
        & $n_\mathrm{cy}$ & $T_\mathrm{EF}$ [s]
        & $t_\mathrm{EF}$ [s/cycle]
        & Rel. \\
        \hline
        Carroll
        & -- & --
        & 17296 & 14325
        & 0.828
        & 1.000 \\

        $\mathrm{PANN}_{\mathcal{SP}}\:(5)$
        & SoftPlus & $1 \times 5$
        & 17460 & 14985
        & 0.858
        & 1.036 \\

        $\mathrm{PANN}_{\mathcal{SP}}\:(5-5)$
        & SoftPlus & $2 \times 5$
        & 17527 & 16785
        & 0.958
        & 1.156 \\

        $\mathrm{PANN}_{\mathcal{SQP}}\:(5)$
        & SQuarePlus & $1 \times 5$
        & 17494 & 14745
        & 0.843
        & 1.018 \\

        $\mathrm{PANN}_{\mathcal{SQP}}\:(5-5)$
        & SQuarePlus & $2 \times 5$
        & 17385 & 15520
        & 0.893
        & 1.078 \\
        \hline
    \end{tabular}
\end{table}

The Carroll model is used as the reference for the runtime comparison. For the shallow PANN architecture with one hidden layer, both neural-network formulations remain computationally close to the classical constitutive model. The normalized runtime ratio increases from $1.000$ for the Carroll model to $1.036$ for the SoftPlus-based PANN and to $1.018$ for the SQuarePlus-based PANN. Thus, the additional element-force runtime is only $3.6\,\%$ for the SoftPlus formulation and $1.8\,\%$ for the SQuarePlus formulation.

For the ``deeper`` PANN architecture with two hidden layers, the additional computational effort becomes more pronounced. The normalized runtime ratio increases to $1.156$ for the SoftPlus-based formulation and to $1.078$ for the SQuarePlus-based formulation. This indicates that the network architecture has a strong influence on the computational cost. At the same time, replacing SoftPlus by SQuarePlus consistently reduces the additional runtime. Relative to the Carroll reference, the overhead decreases from $3.6\,\%$ to $1.8\,\%$ for the shallow architecture and from $15.6\,\%$ to $7.8\,\%$ for the ``deeper`` architecture. In both cases, the additional cost introduced by the PANN formulation is therefore reduced by approximately one half.

The observed improvement in the complete finite element simulation is smaller than the speedup obtained in the isolated activation function benchmark. This is expected, since the element-force calculation contains not only the activation evaluation, but also tensor operations, invariant calculations, stress reconstruction, memory access, and further element-level operations. These contributions partially mask the performance gain of a cheaper activation function. In addition, the relative performance of SoftPlus and SQuarePlus depends on compiler optimizations, architecture-specific code generation, vectorization, and the solver setup. The quantitative runtime ratios should therefore not be interpreted as universal values.

Nevertheless, the simulations consistently show that the SQuarePlus activation reduces the computational overhead of the PANN material routines without changing the overall simulation workflow. The effect becomes more relevant for larger network architectures, where the number of activation evaluations increases. The results therefore indicate that compact PANN architectures can be integrated into explicit finite element simulations with manageable additional cost, and that SQuarePlus is a promising activation function for efficient PANN-based constitutive models.

\section{Conclusion and outlook}
\label{sec:Conclusion}

This work presented a workflow for integrating PANN constitutive models into an explicit finite element simulation environment. Starting from a pretrained PANN, a standalone Fortran user material routine was generated automatically and used in \textsc{Radioss} without relying on external machine-learning libraries during the simulation. The generated routine contains the trained network parameters, the activation functions, the invariant-based stress evaluation, and the volumetric penalty contribution required for the nearly incompressible hyperelastic formulation.

Computational examples demonstrate that the generated PANN material routines can be applied to a complex explicit impact simulation involving large deformations, contact, and a large number of constitutive evaluations. Compared with the Carroll reference model, the PANN formulations produced qualitatively similar stress distributions and did not lead to numerical instabilities. This confirms the feasibility of using automatically generated PANN routines as user-defined hyperelastic material models in an explicit finite element workflow.

The runtime evaluation showed that the additional computational cost of the PANN formulation can remain small if compact network architectures are used. For the shallow PANN with one hidden layer and five neurons, the normalized element-force runtime increased by only $3.6\,\%$ for the SoftPlus activation and by $1.8\,\%$ for the SQuarePlus activation compared with the Carroll model. For the ``deeper`` network with two hidden layers, the overhead increased to $15.6\,\%$ for SoftPlus and $7.8\,\%$ for SQuarePlus. Thus, replacing SoftPlus by SQuarePlus reduced the additional runtime of the PANN material routine by approximately one half in the investigated simulations. Overall, the PANNs are consistently slower than the Carroll model. However, this additional computational cost should be interpreted in relation to their substantially larger number of trainable parameters. The PANN with one hidden layer contains $20$ trainable parameters, whereas the PANN with two hidden layers contains $50$ trainable parameters, compared with only $3$ material parameters in the Carroll model. This increased parameterization provides the PANNs with greater flexibility to represent complex material behavior than the Carroll model.

The quantitative runtime ratios should not be interpreted as universal values, since the relative performance depends on the compiler settings, processor architecture, vectorization strategy, and solver environment. Nevertheless, the results show that the choice of activation function can have a measurable influence on the overall simulation time, particularly for larger network architectures. The SQuarePlus activation therefore represents a promising alternative for efficient PANN-based constitutive models in explicit finite element simulations.

The presented implementation concept is not restricted to the simple hyperelastic PANN architecture considered in this work. Since the constitutive response is evaluated inside the user routine, the framework can in principle be extended to material models with internal variables and history-dependent evolution equations. Future work should therefore focus on generalizing the code-generation procedure to broader PANN architectures and more advanced constitutive formulations, including dissipative and path-dependent material behavior. In addition, the computation of consistent tangent stiffness matrices should be incorporated into the code-generation framework. This would enable the same automatically generated PANN material routines to be used not only in explicit simulations, but also in implicit finite element analyses. Further applications to different material classes and experimentally calibrated datasets would also be useful to assess the robustness and practical applicability of the proposed workflow. Overall, the presented results indicate that physics-augmented neural-network constitutive models can be incorporated into established explicit finite element workflows with manageable computational overhead. The proposed code-generation approach therefore provides a practical step toward the use of machine-learning-based constitutive models in engineering simulations.

\vspace{1cm}

\textbf{Acknowledgment}\\
The authors gratefully acknowledge the contribution of Mehrshad Ilchi, who implemented preliminary ideas for incorporating a PANN architecture into the commercial and open-source explicit finite element codes \textsc{Radioss} and \textsc{OpenRadioss}.

\textbf{Data availability}\\
The Fortran user material routines, stored PANN parameter files, activation-function benchmark code, and the Python code-generation script are available in the public GitHub repository: \url{https://github.com/OVGU-CoMe/PANN_Radioss}. The larger simulation data, including the surface meshes, the \textsc{Radioss} simulation input files for the three mesh resolutions, and the animation shown in Appendix~C, are archived on Zenodo: \url{https://doi.org/10.5281/zenodo.20763660}.

\textbf{Declaration on the use of AI tools}\\
This manuscript was prepared by the authors. ChatGPT and DeepL were used solely for language refinement and stylistic editing.
\newpage
\appendix
\renewcommand{\thefigure}{\Alph{section}-\arabic{figure}}
\counterwithin*{figure}{section}
\renewcommand{\thelstlisting}{\Alph{section}-\arabic{lstlisting}}
\counterwithin*{lstlisting}{section}
\section{Carroll material routine}
\label{App:Subroutine_Carroll}
\revIV{This appendix provides the Fortran user subroutine for the Carroll material model used in \textsc{Radioss}.}
\begin{lstlisting}[
    style=fortranstyle,
    caption={Material user subroutine for the Carroll model},
    label={lst:subrout_carroll},
    breaklines=true,
]
module ca_model_001
  implicit none
  integer, parameter :: dp = selected_real_kind(15, 307)

  real(dp), parameter :: Ca = 1.5e-1_dp
  real(dp), parameter :: Cb = 3.1e-7_dp
  real(dp), parameter :: Cc = 9.5e-2_dp
 
  ! Output shear modulus estimate
  real(dp), parameter :: G = 2.0_dp*Ca + 216.0_dp*Cb + Cc/sqrt(3.0_dp)
  
  ! Volumetric stiffness (to achive the same Speed of Sound as the PANN)
  real(dp), parameter :: K  = 3.00000000000000000e+02_dp

contains

  subroutine pann_eval_from_F(F, sig, K_out, G_out)
    implicit none
    real(dp), intent(in)  :: F(9)
    real(dp), intent(out) :: sig(6), K_out, G_out

    real(dp) :: identity(6)
	real(dp) :: b(6), b2(6)
    real(dp) :: I1, I2, I1b, I2b, J, Jm23, Jm43
    real(dp) :: dI1(6), dI1b(6), dI2b(6)
	real(dp) :: dpsi_dI1b, dpsi_dI2b

    K_out = K
    G_out = G

    ! b in Voigt: [b11,b22,b33,b12,b23,b13]
    b(1) = F(1)*F(1) + F(4)*F(4) + F(6)*F(6)
    b(2) = F(7)*F(7) + F(2)*F(2) + F(5)*F(5)
    b(3) = F(9)*F(9) + F(8)*F(8) + F(3)*F(3)
    b(4) = F(1)*F(7) + F(4)*F(2) + F(5)*F(6)
    b(5) = F(7)*F(9) + F(8)*F(2) + F(3)*F(5)
    b(6) = F(1)*F(9) + F(4)*F(8) + F(3)*F(6)

    ! b^2 in Voigt
    b2(1) = b(1)*b(1) + b(4)*b(4) + b(6)*b(6)
    b2(2) = b(4)*b(4) + b(2)*b(2) + b(5)*b(5)
    b2(3) = b(6)*b(6) + b(5)*b(5) + b(3)*b(3)
    b2(4) = b(1)*b(4) + b(2)*b(4) + b(5)*b(6)
    b2(5) = b(4)*b(6) + b(5)*b(2) + b(3)*b(5)
    b2(6) = b(1)*b(6) + b(4)*b(5) + b(3)*b(6)

    identity = [1.0_dp, 1.0_dp, 1.0_dp, 0.0_dp, 0.0_dp, 0.0_dp]
	
    I1  = b(1) + b(2) + b(3)
    dI1 = identity

    I2  = b(1)*b(2) + b(1)*b(3) + b(2)*b(3) - b(4)**2 - b(5)**2 - b(6)**2
    
    J = F(1)*F(2)*F(3) + F(4)*F(5)*F(9) + F(6)*F(7)*F(8) - F(9)*F(2)*F(6) - F(8)*F(5)*F(1) - F(3)*F(7)*F(4)
    Jm23 = J**(-2.0_dp/3.0_dp)
    Jm43 = Jm23*Jm23
	
    I1b  = Jm23 * I1
    dI1b = Jm23*(b-1.0_dp/3.0_dp*I1*dI1)

    I2b  = Jm43 * I2
    dI2b = Jm43*(I1*b - b2 - (2.0_dp/3.0_dp)*I2*dI1)

    dpsi_dI1b = Ca + 4.0_dp*Cb*I1b**3
    dpsi_dI2b = Cc/(2.0_dp*sqrt(I2b))
	
    ! Isochoric Carroll stress contribution
    sig = 2.0_dp * (dpsi_dI1b*dI1b + dpsi_dI2b*dI2b)

    ! Volume term for psi_vol = 0.5*K*(J-1)^2
    sig = sig/J + K*(J-1.0_dp)*dI1

  end subroutine pann_eval_from_F

end module ca_model_001
\end{lstlisting}

\newpage
\section{Benchmark SoftPlus vs. SQuarePlus}
\label{App:Benchmark_SoftSquare}
This appendix gives the computational comparison between the SoftPlus and the SQuarePlus function.

\begin{lstlisting}[
    style=fortranstyle,
    caption={Benchmark SoftPlus vs. SQuarePlus},
    label={lst:benchm_softpl_square},
    breaklines=true,
]
! Benchmark for SoftPlus and SQuarePlus activation functions
!
! The program compares the isolated evaluation cost of two activation functions
! and their first derivatives. 
! -----------------------------------------------------------------------------
program bm_s_sq
    use iso_fortran_env, only: int64, real64
    implicit none

    integer, parameter :: dp = real64
    integer(int64), parameter :: n = 100000000_int64   ! Number of evaluations per benchmark loop

    integer(int64) :: i
    integer(int64) :: c_start, c_end, rate
    real(dp) :: x, y, dy
    real(dp) :: sum_soft, sum_square
    real(dp) :: time_soft, time_square
    real(dp) :: speedup, saving_percent
    real(dp) :: time_base
    real(dp) :: time_soft_total, time_square_total
    real(dp) :: time_soft_net, time_square_net
    real(dp) :: sum_base
    
	! Determine the resolution of the system clock used for all measurements.
    call system_clock(count_rate=rate)

    print *, "Number of evaluations per function:", n
    print *, "System clock rate:", rate, "ticks/s"
    print *, ""

    ! Baseline measurement
    ! ------------------------------------------------------------
    ! This loop performs the same x-value generation and accumulation as the
    ! activation benchmarks, but without calling an activation function. The
    ! measured runtime is subtracted later to estimate the net function cost.
    sum_base = 0.0_dp

    call system_clock(c_start)

    do i = 1_int64, n
	    ! x is sampled uniformly from approximately -20 to +20.
        x = -20.0_dp + 40.0_dp * real(i - 1_int64, dp) / real(n - 1_int64, dp)

        ! Dummy operations matching the accumulation pattern in the benchmark
        ! loops. This keeps the baseline comparable to the two function tests.
        y  = x
        dy = 1.0_dp

        sum_base = sum_base + y + dy
    end do

    call system_clock(c_end)

    time_base = real(c_end - c_start, dp) / real(rate, dp)

    ! SoftPlus benchmark
    ! ------------------------------------------------------------
    ! The activation value and its derivative are evaluated for the same x
    ! values as in the baseline loop.
    sum_soft = 0.0_dp

    call system_clock(c_start)

    do i = 1_int64, n
        x = -20.0_dp + 40.0_dp * real(i - 1_int64, dp) / real(n - 1_int64, dp)

        call softplus_both(x, y, dy)

        ! Accumulation prevents the compiler from treating y and dy as unused.
        sum_soft = sum_soft + y + dy
    end do

    call system_clock(c_end)

    time_soft_total = real(c_end - c_start, dp) / real(rate, dp)

    ! SQuarePlus benchmark
    ! ------------------------------------------------------------
    ! The same x values and accumulation pattern are used to make the comparison
    ! with the SoftPlus benchmark as direct as possible.
    sum_square = 0.0_dp

    call system_clock(c_start)

    do i = 1_int64, n
        x = -20.0_dp + 40.0_dp * real(i - 1_int64, dp) / real(n - 1_int64, dp)

        call squareplus_both(x, y, dy)

        ! Accumulation prevents the compiler from treating y and dy as unused.
        sum_square = sum_square + y + dy
    end do

    call system_clock(c_end)

    time_square_total = real(c_end - c_start, dp) / real(rate, dp)

    ! Remove the baseline runtime to estimate the isolated activation-function
    ! evaluation time. These net values are used for the reported speedup.
    time_soft_net   = time_soft_total   - time_base
    time_square_net = time_square_total - time_base

    print *, "Control sum baseline :", sum_base
    print *, "Control sum softplus  :", sum_soft
    print *, "Control sum squareplus:", sum_square
    print *, ""

    print *, "Baseline time      :", time_base, "seconds"
    print *, ""

    print *, "Softplus gross time    :", time_soft_total, "seconds"
    print *, "Squareplus gross time  :", time_square_total, "seconds"

    ! Compute the net speedup only if both baseline-corrected timings are valid.
    if (time_soft_net > 0.0_dp .and. time_square_net > 0.0_dp) then
        speedup = time_soft_net / time_square_net
        saving_percent = (1.0_dp - time_square_net / time_soft_net) * 100.0_dp

        print *, "Net speedup of squareplus compared to softplus:", speedup
        print *, "Net time saving with squareplus:", saving_percent, "%"
    else
        print *, "Net times are too small or invalid."
        print *, "Increase n or perform multiple repetitions."
    end if

contains

    ! SoftPlus activation and first derivative
    ! ------------------------------------------------------------
    pure elemental subroutine softplus_both(x, y, dy)
        real(dp), intent(in)  :: x
        real(dp), intent(out) :: y, dy
        real(dp) :: ax, t

        ax = abs(x)
        t  = exp(-ax)
        y  = max(x, 0.0_dp) + log(1.0_dp + t)
        dy = merge( 1.0_dp / (1.0_dp + t), &
                    t      / (1.0_dp + t), &
                    x >= 0.0_dp )
    end subroutine softplus_both

    ! SQuarePlus activation and first derivative
    ! ------------------------------------------------------------
    pure elemental subroutine squareplus_both(x, y, dy)
        real(dp), intent(in)  :: x
        real(dp), intent(out) :: y, dy
        real(dp), parameter   :: b = 2.0_dp
        real(dp) :: s

        s  = sqrt(x*x + b)
        y  = 0.5_dp * (x + s)
        dy = 0.5_dp * (1.0_dp + x / s)
    end subroutine squareplus_both

end program bm_s_sq
\end{lstlisting}

\newpage
\section{Simulation sequence}
\label{App:Sequence}
\begin{figure}[htb]
    \centering
    \def\svgwidth{0.80\linewidth}
    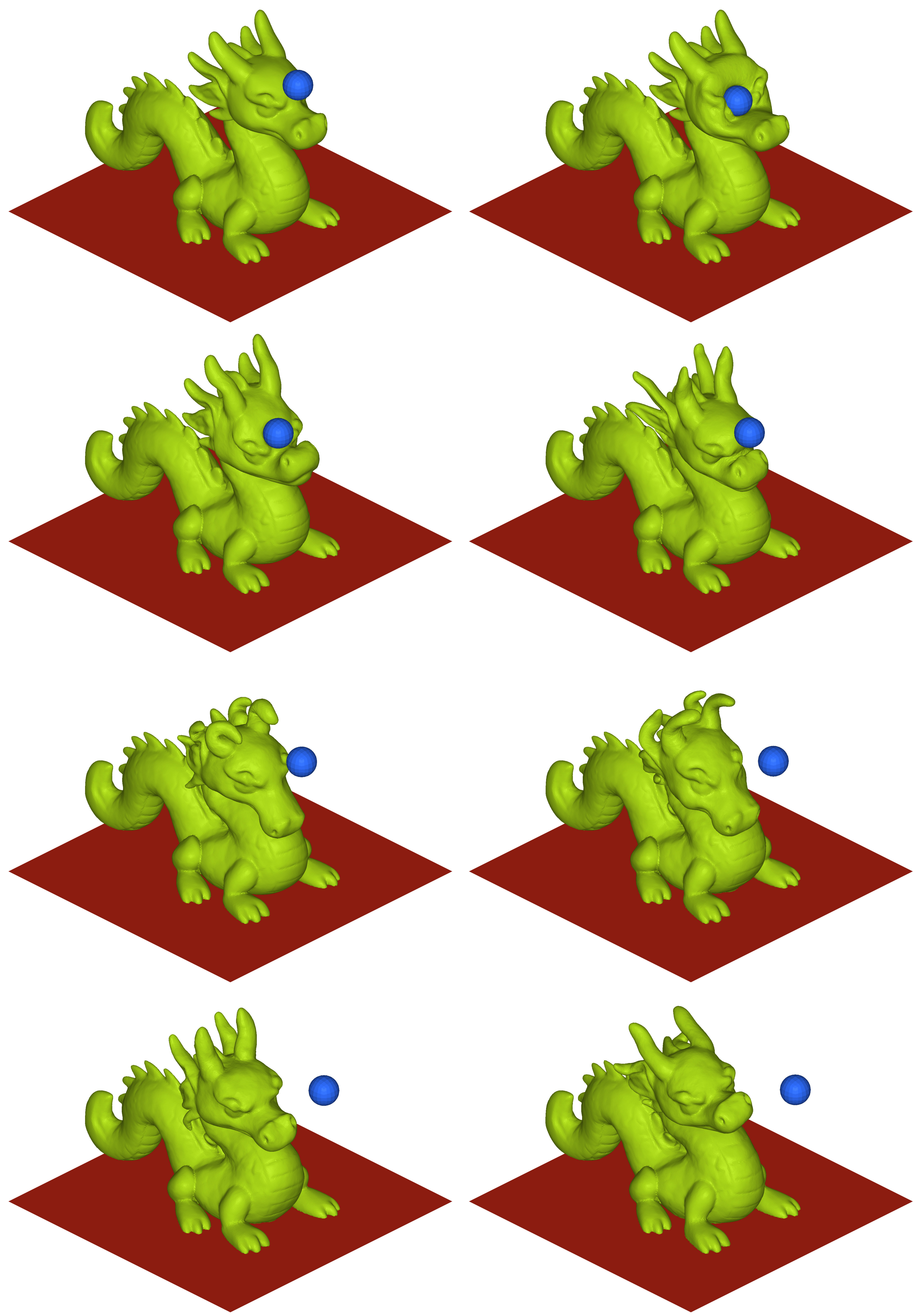
    \caption{Simulation sequence of the dragon benchmark using Mesh~03 and the Carroll material model.}
    \label{fig:Sequence}
\end{figure}

The benchmark simulation is illustrated in Fig.~\ref{fig:Sequence} for Mesh~3 using the Carroll material model.

\cleardoublepage
\bibliographystyle{unsrt}
\phantomsection
\addcontentsline{toc}{section}{References}
\bibliography{./bib/references}
\end{document}